\documentclass[10.5pt]{article}

\def\sqr#1#2{{\vcenter{\vbox{\hrule height.#2pt
              \hbox{\vrule width.#2pt height#1pt \kern#1pt \vrule width.#2pt}
              \hrule height.#2pt}}}}
\def\signed #1{{\unskip\nobreak\hfil\penalty50
              \hskip2em\hbox{}\nobreak\hfil#1
              \parfillskip=0pt \finalhyphendemerits=0 \par}}
\def\endpf{\signed {$\sqr69$}}

\def\3n{\negthinspace \negthinspace \negthinspace }
\def\2n{\negthinspace \negthinspace }
\def\1n{\negthinspace }

\def\dbE{{\mathop{\rm l\negthinspace E}}}
\def\dbF{{\mathop{\rm l\negthinspace F}}}

\def\dbP{{\mathop{\rm l\negthinspace P}}}
\def\dbR{{\mathop{\rm l\negthinspace R}}}
\def\={\buildrel \triangle \over =}

\def\ds{\displaystyle}

\def\ns{\noalign{\ss}}

\def\a{\alpha}

\def\d{\delta}
\def\e{\varepsilon}
\def\z{\zeta}

\def\t{\tau}
\def\f{\varphi}

\def\o{\omega}

\def\G{\Gamma}
\def\D{\Delta}
\def\Th{\Theta}

\def\O{\Omega}

\def\cA{{\cal A}}
\def\cB{{\cal B}}

\def\cF{{\cal F}}
\def\cG{{\cal G}}

\def\cQ{{\cal Q}}
\def\cR{{\cal R}}
\def\cS{{\cal S}}

\def\cU{{\cal U}}

\def\cX{{\cal X}}

\def\ss{\smallskip}
\def\ms{\medskip}

\def\q{\quad}
\def\qq{\qquad}
\def\hb{\hbox}

\def\lan{\mathop{\langle}}
\def\ran{\mathop{\rangle}}

\def\h{\widehat}
\def\wt{\widetilde}
\def\cd{\cdot}
\def\cds{\cdots}

\def\as{\hbox{\rm a.s.{ }}}

\def\var{\hbox{\rm var$\,$}}

\def\deq{\mathop{\buildrel\D\over=}}

\def\({\Big (}
\def\){\Big )}
\def\[{\Big[}
\def\]{\Big]}
\def\bde{\begin{definition}}
\def\ede{\end{definition}}
\def\be{\begin{equation}}
\def\bel{\begin{equation}\label}
\def\ee{\end{equation}}
\def\bt{\begin{theorem}}
\def\et{\end{theorem}}
\def\bc{\begin{corollary}}
\def\ec{\end{corollary}}
\def\bl{\begin{lemma}}
\def\el{\end{lemma}}
\def\bp{\begin{proposition}}
\def\ep{\end{proposition}}
\def\bas{\begin{assumption}}
\def\eas{\end{assumption}}
\def\br{\begin{remark}}
\def\er{\end{remark}}
\def\ba{\begin{array}}
\def\ea{\end{array}}
\def\ed{\end{document}}

\def\square#1{\vbox{\hrule\hbox{\vrule height#1%
     \kern#1\vrule}\hrule}}
\def\rectangle#1#2{\vbox{\hrule\hbox{\vrule height#1%
     \kern#2\vrule}\hrule}}

\font\tenbb=msbm10 \font\sevenbb=msbm7 \font\fivebb=msbm5

\newfam\bbfam
\scriptscriptfont\bbfam=\fivebb \textfont\bbfam=\tenbb
\scriptfont\bbfam=\sevenbb

\newtheorem{lemma}{Lemma}[section]
\newtheorem{remark}{Remark}[section]

\newtheorem{theorem}{Theorem}[section]
\newtheorem{corollary}{Corollary}[section]

\newtheorem{definition}{Definition}[section]
\newtheorem{proposition}{Proposition}[section]
\newtheorem{assumption}{Assumption}[section]

\makeatletter
   
   \@addtoreset{equation}{section}
\makeatother

\begin{document}
\title{\bf A Linear-Quadratic
Optimal Control Problem\\ for Mean-Field Stochastic Differential
Equations\footnote{This work is supported in part by NSF Grant
DMS-1007514.}}

\author{Jiongmin Yong\\ Department of
Mathematics, University of Central Florida, Orlando, FL 32816, USA.}

\maketitle

\begin{abstract} A Linear-quadratic optimal control problem is considered for mean-field stochastic
differential equations with deterministic coefficients. By a
variational method, the optimality system is derived, which turns
out to be a linear mean-field forward-backward stochastic
differential equation. Using a decoupling technique, two Riccati
differential equations are obtained, which are uniquely solvable
under certain conditions. Then a feedback representation is obtained
for the optimal control.

\end{abstract}

\ms


\ms


\section{Introduction.}

Let $(\O,\cF,\dbP,\dbF)$ be a complete filtered probability space,
on which a one-dimensional standard Brownian motion $W(\cd)$ is
defined with $\dbF\equiv\{\cF_t\}_{t\ge0}$ being its natural
filtration augmented by all the $\dbP$-null sets. Consider the
following controlled linear stochastic differential equation (SDE,
for short):
\bel{1.1}\left\{\3n\ba{ll}
\ns\ds dX(s)=\Big\{A(s)X(s)+B(s)u(s)+\h A(s)\dbE[X(s)]+\h
B(s)\dbE[u(s)]\Big\}ds\\
\ns\ds\qq\qq+\Big\{A_1(s)X(s)\1n+\1n B_1(s)u(s)\1n+\1n\h
A_1(s)\dbE[X(s)]\1n+\1n\h
B_1(s)\dbE[u(s)]\Big\}dW(s),\q s\in[0,T],\\
\ns\ds X(0)=x,\ea\right.\ee
where $A(\cd),B(\cd),\h A(\cd),\h B(\cd),A_1(\cd),B_1(\cd),\h
A_1(\cd),\h B_1(\cd)$ are given deterministic matrix valued
functions. In the above, $X(\cd)$, valued in $\dbR^n$, is the {\it
state process}, and $u(\cd)$, valued in $\dbR^m$, is the {\it
control process}.

\ms

We note that $\dbE[X(\cd)]$ and $\dbE[u(\cd)]$ appear in the state
equation. Such an equation is referred to as a linear controlled
mean-field (forward) SDE (MF-FSDE, for short). MF-FSDEs can be used
to describe particle systems at mesoscopic level, which is of great
importance in applications. Historically, the later-called
McKean--Vlasov SDE, a kind of MF-FSDE, was suggested by Kac
\cite{Kac 1956} in 1956 as a stochastic toy model for the Vlasov
kinetic equation of plasma and the study of which was initiated by
McKean \cite{McKean 1966} in 1966. Since then, many authors made
contributions on McKean--Vlasov type SDEs and applications, see, for
examples, Dawson \cite{Dawson 1983}, Dawson--G\"artner
\cite{Dawson-Gartner 1987}, G\'artner \cite{Gartner 1988}, Scheutzow
\cite{Scheutzow 1987}, Graham \cite{Graham 1992}, Chan \cite{Chan
1994}, Chiang \cite{Chiang 1994}, Ahmed--Ding \cite{Ahmed-Ding
1995}, to mention a few. In recent years, related topics and
problems have attracted more and more researchers' attentions, see,
for examples, Veretennikov \cite{Veretennikov 2003},
Huang--Malham\'e--Caines \cite{Huang-Malhame-Caines 2006},
Mahmudov--McKibben \cite{Mahmudov-Mckibben 2007},
Buckdahn-Djehiche-Li-Peng \cite{Buckdahn-Djehiche-Li-Peng 2009},
Buckdahn-Li-Peng \cite{Buckdahn-Li-Peng 2009}, Borkar--Kumar
\cite{Borkar-Kumar 2010}, Crisan--Xiong \cite{Crisan-Xiong 2010},
Kotelenez--Kurtz \cite{Kotelenez-Kurtz 2010}, Kloeden--Lorenz
\cite{Kloeden-Lorenz 2010}, and so on. More interestingly, control
problems of McKean--Vlasov equation or MF-FSDEs were investigated by
Ahmed--Ding \cite{Ahmed-Ding 2001}, Ahmed \cite{Ahmed 2007},
Buckdahn--Djehiche--Li \cite{Buckdahn-Djehiche-Li 2010},
Park--Balasubramaniam--Kang \cite{Park-Balasubramaniam-Kang 2008},
Andersson--Djehiche \cite{Andersson-Djehche 2011},
Meyer-Brandis--Oksendal--Zhou \cite{Meyer-Brandis-Oksandal-Zhou
2011}, and so on. This paper can be regarded as an addition to the
study of optimal control for MF-FSDEs.

\ms

For the state equation (\ref{1.1}), we introduce the following:
$$\cU[0,T]=L^2_\dbF(0,T;\dbR^m)\deq\Big\{u:[0,T]\times\O\to\dbR^m\bigm|\dbE\int_0^T|u(s)|^2ds<\infty\Big\}.$$
Any $u(\cd)\in\cU[0,T]$ is called an {\it admissible control}. Under
mild conditions, one can show that (see below) for any
$(x,u(\cd))\in\dbR^n\times\cU[0,T]$, (\ref{1.1}) admits a unique
solution $X(\cd)=X(\cd\,;x,u(\cd))$. We introduce the following cost
functional:
\bel{1.2}\ba{ll}
\ns\ds J(x;u(\cd))=\dbE\Big\{\int_0^T\[\lan Q(s)X(s),X(s)\ran+\lan\h
Q(s)\dbE[X(s)],\dbE[X(s)]\ran+\lan R(s)u(s),u(s)\ran\\
\ns\ds\qq\qq\qq\qq\q+\lan\h R(s)\dbE[u(s)],\dbE[u(s)]\ran\]ds+\lan
GX(T),X(T)\ran+\lan\h G\dbE[X(T)],\dbE[X(T)]\ran\Big\},\ea\ee
with $Q(\cd),R(\cd),\h Q(\cd),\h R(\cd)$ being suitable
deterministic symmetric matrix-valued functions, and $G,\h G$ being
symmetric matrices. Our optimal control problem can be stated as
follows:

\ms

\bf Problem (MF). \rm For given $x\in\dbR^n$, find a $\bar
u(\cd)\in\cU[0,T]$ such that
\bel{}J(x;\bar u(\cd))=\inf_{u(\cd)\in\cU[0,T]}J(x;u(\cd)).\ee

\ms

Any $\bar u(\cd)\in\cU[0,T]$ satisfying the above is called an {\it
optimal control} and the corresponding state process $\bar
X(\cd)\equiv X(\cd\,;x,\bar u(\cd))$ is called an {\it optimal state
process}; the pair $(\bar X(\cd),\bar u(\cd))$ is called an {\it
optimal pair}.

\ms

From the above-listed literature, one has some motivations for the
inclusion of $\dbE[X(\cd)]$ and $\dbE[u(\cd)]$ in the state
equation. We now briefly explain a motivation of including
$\dbE[X(\cd)]$ and $\dbE[u(\cd)]$ in the cost functional. We recall
that for a classical LQ problem with the state equation
\bel{1.4}\left\{\3n\ba{ll}
\ns\ds
dX(s)=\[A(s)X(s)+B(s)u(s)\]ds+\[A_1(s)X(s)+B_1(s)u(s)\]dW(s),\q s\in[0,T],\\
\ns\ds X(0)=x,\ea\right.\ee
one has the following cost functional
\bel{1.5}J_0(x;u(\cd))=\dbE\Big\{\int_0^T\[\lan
Q_0(s)X(s),X(s)\ran+\lan R_0(s)u(s),u(s)\ran\]ds+\lan
G_0X(T),X(T)\ran\Big\}.\ee
For such a corresponding optimal control problem, it is natural to
hope that the optimal state process and/or control process are not
too sensitive with respect to the possible variation of the random
events. One way to achieve this is try to make the variations
$\var[X(\cd)]$ and $\var[u(\cd)]$ small. Therefore, one could
include the $\var[X(\cd)]$ and $\var[u(\cd)]$ in the cost
functional. Consequently, one might want to replace (\ref{1.5}) by
the following:
\bel{1.5}\ba{ll}
\ns\ds\h J_0(x;u(\cd))=\dbE\Big\{\int_0^T\[\lan
Q_0(s)X(s),X(s)\ran+q(s)\var[X(s)]+\lan
R_0(s)u(s),u(s)\ran+\rho(s)\var[u(s)]\]ds\\
\ns\ds\qq\qq\qq\qq\qq+\lan G_0X(T),X(T)\ran+g\var[X(T)]\Big\},\ea\ee
for some (positive) weighting factors $q(\cd)$, $\rho(\cd)$, and
$g$. Since
$$\var[X(s)]=\dbE|X(s)|^2-\(\dbE[X(s)]\)^2,$$
and similar things hold for $\var[X(T)]$ and $\var[u(s)]$, we see
that
$$\ba{ll}
\ns\ds\h J_0(x;u(\cd))=\dbE\Big\{\int_0^T\[\lan
[Q_0(s)+q(s)I]X(s),X(s)\ran-q(s)\(\dbE[X(s)]\)^2+\lan
[R_0(s)+\rho(s)I]u(s),u(s)\ran\\
\ns\ds\qq\qq\qq\qq\qq-\rho(s)\(\dbE[u(s)]\)^2\]ds+\lan
[G_0+gI]X(T),X(T)\ran-g\(\dbE[X(T)]\)^2\Big\}.\ea$$
Clearly, the above is a special case of (\ref{1.2}) with
$$\ba{ll}
\ns\ds Q(\cd)=Q_0(\cd)+q(\cd)I,\q R(\cd)=R_0(\cd)+\rho(\cd)I,\q G=G_0+gI,\\
\ns\ds\h Q(\cd)=-q(\cd)I,\q\h R(\cd)=-\rho(\cd)I,\q\h G=-gI.\ea$$
Note that in the above case, $\h Q(\cd)$, $\h R(\cd)$, and $\h G$
are not positive semi-definite. Because of this, we do not assume
the positive semi-definiteness for $\h Q(\cd)$, $\h R(\cd)$, and $\h
G$.

\ms

The purpose of this paper is to study Problem (MF). We will begin
with the well-posedness of the state equation and the solvability of
Problem (MF) in Section 2. Then, in Section 3, we will establish
necessary and sufficient conditions for optimal pairs. Naturally, a
linear backward stochastic differential equation of mean-filed type
(MF-BSDE, for short) will be derived. Consequently, the optimality
system turns out to be a coupled mean-field type forward-backward
stochastic differential equation (MF-FBSDE, for short). Inspired by
the invariant imbedding \cite{Bellman 1960}, and the Four-Step
Scheme \cite{MY 1999}, we derive two Riccati differential equations
in Section 4, so that the optimal control is represented as a state
feedback form. Well-posedness of these Riccati equations will be
established. We also present a direct verification of optimality for
the state feedback control by means of completing squares. In
Section 5, we will look at a modified LQ problem which is one
motivation of the current paper.

\section{Preliminaries.}

\ms

In this section, we make some preliminaries. First of all, for any
Euclidean space $H=\dbR^n,\dbR^{n\times m},\cS^n$ (with $\cS^n$
being the set of all $(n\times n)$ symmetric matrices), we let
$L^p(0,t;H)$ be the set of all $H$-valued functions that are
$L^p$-integrable, $p\in[1,\infty]$. Next, we introduce the following
spaces:
$$\ba{ll}
\ns\ds\cX_t\equiv
L^2_{\cF_t}(\O;\dbR^n)=\Big\{\xi:\O\to\dbR^n\bigm|\xi\hb{ is
$\cF_t$-measurable, }\dbE|\xi|^2<\infty\Big\},\\
\ns\ds\cU_t\equiv
L^2_{\cF_t}(\O;\dbR^m)=\Big\{\eta:\O\to\dbR^m\bigm|\eta\hb{ is
$\cF_t$-measurable, }\dbE|\eta|^2<\infty\Big\},\\
\ns\ds
L^2_\dbF(0,t;\dbR^n)=\Big\{X:[0,t]\times\O\to\dbR^n\bigm|X(\cd)\hb{
is $\dbF$-adapted,
}\dbE\int_0^t|X(s)|^2ds<\infty\Big\},\\
\ns\ds\cX[0,t]\equiv
C_\dbF([0,t];\dbR^n)=\Big\{X:[0,t]\to\cX_t\bigm|X(\cd)\hb{ is
$\dbF$-adapted,}\\
\ns\ds\qq\qq\qq\qq\qq\qq\q\hb{ $s\mapsto X(s)$ is continuous, }
\sup_{s\in[0,t]}\dbE|X(s)|^2<\infty\Big\},\\
\ns\ds\h\cX[0,t]\equiv
L^2_\dbF(\O;C([0,t];\dbR^n))=\Big\{X:[0,t]\times\O\to\dbR^n\bigm|X(\cd)\hb{
is $\dbF$-adapted, }\\
\ns\ds\qq\qq\qq\qq\qq\qq\qq\q X(\cd)\hb{ has continuous paths,
}\dbE\[\sup_{s\in[0,t]}|X(s)|^2\]<\infty\Big\}.\ea$$
Note that in the definition of $\cX[0,t]$, the continuity of
$s\mapsto X(s)$ means that as a map from $[0,t]$ to $\cX_t$, it is
continuous. Whereas, in the definition of $\h\cX[0,t]$, $X(\cd)$ has
continuous paths means that for almost sure $\o\in\O$, $s\mapsto
X(s,\o)$ is continuous. It is known that
$$\ba{ll}
\ns\ds\h\cX[0,t]\subseteq\cX[0,t]\subseteq L^2_\dbF(0,t;\dbR^n),\qq
\h\cX[0,t]\ne\cX[0,t]\ne L^2_\dbF(0,t;\dbR^n).\ea$$
We now introduce the following assumptions for the coefficients of
the state equation.

\ms

{\bf(H1)} The following hold:
\bel{H1}\left\{\ba{ll}
\ns\ds A(\cd),\h A(\cd),A_1(\cd),\h A_1(\cd)\in
L^\infty(0,T;\dbR^{n\times n}),\\
\ns\ds B(\cd),\h B(\cd),B_1(\cd),\h B_1(\cd)\in
L^\infty(0,T;\dbR^{n\times m}).\ea\right.\ee

\ms

{\bf(H1)$'$} The following hold:
\bel{H1'}\left\{\ba{ll}
\ns\ds A(\cd),\h A(\cd)\in L^2(0,T;\dbR^{n\times n}),\qq A_1(\cd),\h
A_1(\cd)\in L^\infty(0,T;\dbR^{n\times n}),\\
\ns\ds B(\cd),\h B(\cd)\in L^2(0,T;\dbR^{n\times m}),\qq B_1(\cd),\h
B_1(\cd)\in L^\infty(0,T;\dbR^{n\times m}).\ea\right.\ee

\ms

{\bf(H1)$''$} The following hold:
\bel{H1''}\left\{\ba{ll}
\ns\ds A(\cd),\h A(\cd)\in L^1(0,T;\dbR^{n\times n}),\qq A_1(\cd),\h
A_1(\cd)\in L^2(0,T;\dbR^{n\times n}),\\
\ns\ds B(\cd),\h B(\cd)\in L^2(0,T;\dbR^{n\times m}),\qq B_1(\cd),\h
B_1(\cd)\in L^\infty(0,T;\dbR^{n\times m}).\ea\right.\ee

\ms

Clearly, (H1) implies (H1)$'$ which further implies (H1)$''$.
Namely, (H1)$''$ is the weakest assumption among the above three.
Whereas, (H1) is the most common assumption. For the weighting
matrices in the cost functional, we introduce the following
assumption.

\ms

{\bf(H2)} The following hold:
\bel{H2}Q(\cd),\h Q(\cd)\in L^\infty(0,T;\cS^n),\q R(\cd),\h
R(\cd)\in L^\infty(0,T;\cS^m),\q G,\h G\in\cS^n.\ee

\ms

{\bf(H2)$'$} In addition to (H2), the following holds:
\bel{}\left\{\ba{ll}
\ns\ds Q(s),\;Q(s)+\h Q(s)\ge0,\qq R(s),\;R(s)+\h R(s)\ge0,\qq
s\in[0,T],\\
\ns\ds G,\;G+\h G\ge0.\ea\right.\ee

\ms

{\bf(H2)$''$} In addition to (H2), the following holds: For some
$\d>0$,
\bel{}\left\{\ba{ll}
\ns\ds Q(s),\;Q(s)+\h Q(s)\ge0,\qq R(s),\;R(s)+\h R(s)\ge\d I,\qq
s\in[0,T],\\
\ns\ds G,\;G+\h G\ge0.\ea\right.\ee

From (\ref{1.4}), we see that $\h Q(\cd)$, $\h R(\cd)$, and $\h G$
are not necessarily positive semi-definite. Therefore, in (H2), we
do not mention positive-definiteness of the involved matrices and
matrix-valued functions.

\ms

Now, for any $X(\cd)\in L^2_\dbF(0,T;\dbR^n)$ and any
$u(\cd)\in\cU[0,T]$, we define
\bel{AB}\left\{\ba{ll}
\ns\ds[\cA X(\cd)](t)=\int_0^t\(A(s)X(s)+\h
A(s)\dbE[X(s)]\)ds\\
\ns\ds\qq\qq\qq\qq+\int_0^t\(A_1(s)X(s)+\h
A_1(s)\dbE[X(s)]\)dW(s),\qq t\in[0,T],\\
\ns\ds[\cB u(\cd)](t)=\int_0^t\(B(s)u(s)+\h
B(s)\dbE[u(s)]\)ds\\
\ns\ds\qq\qq\qq\qq+\int_0^t\(B_1(s)u(s)+\h
B_1(s)\dbE[u(s)]\)dW(s),\qq t\in[0,T].\ea\right.\ee
The following result is concerned with operators $\cA$ and $\cB$.

\ms

\bf Lemma 2.1. \sl The following estimates hold as long as the
involved norms on the right hand sides are meaningful: For any
$t\in[0,T]$,
\bel{estimate1}\ba{ll}
\ns\ds\|\cA X(\cd)\|^2_{\h\cX[0,t]}\le
C\[\|A(\cd)\|_{L^2(0,t;\dbR^{n\times n})}^2+\|\h
A(\cd)\|_{L^2(0,t;\dbR^{n\times
n})}^2\\
\ns\ds\qq\qq\qq\qq+\|A_1(\cd)\|_{L^\infty(0,t;\dbR^{n\times
n})}^2+\|\h A_1(\cd)\|_{L^\infty(0,t;\dbR^{n\times
n})}^2\]\|X(\cd)\|^2_{L^2_\dbF(0,t;\dbR^n)},\ea\ee
\bel{estimate2}\ba{ll}
\ns\ds\|\cA X(\cd)\|_{\h\cX[0,t]}^2\le
C\[\|A(\cd)\|_{L^1(0,t;\dbR^{n\times n})}^2+\|\h
A(\cd)\|_{L^1(0,t;\dbR^{n\times
n})}^2\\
\ns\ds\qq\qq\qq\qq+\|A_1(\cd)\|_{L^2(0,t;\dbR^{n\times n})}^2+\|\h
A_1(\cd)\|_{L^2(0,t;\dbR^{n\times
n})}^2\]\|X(\cd)\|_{\cX[0,t]}^2,\ea\ee
and
\bel{estimate3}\ba{ll}
\ns\ds\|\cB u(\cd)\|^2_{\h\cX[0,t]}\le
C\[\|B(\cd)\|_{L^2(0,t;\dbR^{n\times n})}^2+\|\h
B(\cd)\|_{L^2(0,t;\dbR^{n\times
n})}^2\\
\ns\ds\qq\qq\qq\qq+\|B_1(\cd)\|_{L^\infty(0,t;\dbR^{n\times
n})}^2+\|\h B_1(\cd)\|_{L^\infty(0,t;\dbR^{n\times
n})}^2\]\|u(\cd)\|^2_{\cU[0,t]}.\ea\ee
Hereafter, $C>0$ represents a generic constant which can be
different from line to line.

\ms

\it Proof. \rm For any $t\in(0,T]$, and any $X(\cd)\in
L^2_\dbF(0,t;\dbR^n)$,
$$\ba{ll}
\ns\ds\dbE\(\sup_{s\in[0,t]}|[\cA X(\cd)](s)|\)\le
C\Big\{\dbE\(\int_0^t|A(s)||X(s)|ds\)^2+\(\int_0^t|\h
A(s)||\dbE[X(s)]|ds\)^2\\
\ns\ds\qq+\dbE\int_0^t|A_1(s)|^2|X(s)|^2ds+\int_0^t|\h
A_1(s)|^2|\dbE[X(s)]|^2ds\Big\}\\
\ns\ds\le
C\Big\{\dbE\[\(\int_0^t|A(s)|^2ds\)\(\int_0^t|X(s)|^2ds\)\]+\(\int_0^t|\h
A(s)|^2ds\)\(\int_0^t|\dbE[X(s)]|^2ds\)\\
\ns\ds\qq+\(\sup_{s\in[0,t]}|A_1(s)|\)^2\int_0^t\dbE|X(s)|^2ds+\(\sup_{s\in[0,t]}|\h
A_1(s)|^2\)\int_0^t\dbE|X(s)|^2ds\Big\}\\
\ns\ds\le C\Big\{\(\int_0^t|A(s)|^2ds\)+\(\int_0^t|\h
A(s)|^2ds\)+\sup_{s\in[0,t]}|A_1(s)|^2+\sup_{s\in[0,t]}|\h
A_1(s)|^2\Big\}\int_0^t\dbE|X(s)|^2ds.\ea$$
Thus, estimate (\ref{estimate1}) holds. Next, for any $X(\cd)\in
\cX[0,t]$, making use of Burkholder--Davis--Gundy's inequality
(\cite{Yong-Zhou 1999}), we have
$$\ba{ll}
\ns\ds\dbE\(\sup_{s\in[0,t]}|[\cA X(\cd)](s)|^2\)\le
C\Big\{\dbE\(\int_0^t|A(s)||X(s)|ds\)^2+\(\int_0^t|\h
A(s)||\dbE[X(s)]|ds\)^2\\
\ns\ds\qq\qq\qq\qq\qq+\dbE\int_0^t|A_1(s)|^2|X(s)|^2ds+\int_0^t|\h
A_1(s)|^2|\dbE[X(s)]|^2ds\Big\}\\
\ns\ds\le
C\Big\{\(\int_0^t|A(s)|ds\)\(\int_0^t|A(s)|\dbE|X(s)|^2ds\)+\(\int_0^t|\h
A(s)|ds\)^2\(\sup_{s\in[0,t]}|\dbE[X(s)]|^2\)\\
\ns\ds\qq\qq\qq\qq\qq+\int_0^t|A_1(s)|^2\dbE|X(s)|^2ds+\int_0^t|\h
A_1(s)|^2\dbE|X(s)|^2ds\Big\}\\
\ns\ds\le C\[\(\int_0^t|A(s)|ds\)^2+\(\int_0^t|\h
A(s)|ds\)^2+\int_0^t|A_1(s)|^2ds+\int_0^t|\h
A_1(s)|^2ds\]\(\sup_{s\in[0,t]}\dbE|X(s)|^2\).\ea$$
Hence, (\ref{estimate2}) follows.

\ms

Similar to the proof of (\ref{estimate1}), we can prove
(\ref{estimate3}). \endpf

\ms

The above lemma leads to the following corollary.

\ms

\bf Corollary 2.2. \sl If {\rm(H1)$''$} holds, then
$\cA:\cX[0,T]\to\h\cX[0,T]$ and $\cB:\cU[0,T]\to\h\cX[0,T]$ are
bounded. Further, if {\rm(H1)$'$} holds, then
$\cA:L^2_\dbF(0,T;\dbR^n)\to\h\cX[0,T]$ is also bounded. In
particular, all the above hold if {\rm(H1)} holds.

\rm

\ms

Next, we define
\bel{AB3}\left\{\ba{ll}
\ns\ds I_TX(\cd)=X(T),\\
\ns\ds\cA_TX(\cd)=I_T\cA X(\cd)\equiv[\cA
X(\cd)](T)\\
\ns\ds\qq\q=\int_0^T\(A(s)X(s)+\h
A(s)\dbE[X(s)]\)ds+\int_0^T\(A_1(s)X(s)+\h
A_1(s)\dbE[X(s)]\)dW(s),\\
\ns\ds\cB_Tu(\cd)=I_T\cB u(\cd)\equiv[\cB
u(\cd)](T)\\
\ns\ds\qq\q=\int_0^T\(B(s)u(s)+\h
B(s)\dbE[u(s)]\)ds+\int_0^T\(B_1(s)u(s)+\h
B_1(s)\dbE[u(s)]\)dW(s).\ea\right.\ee
It is easy to see that
$$I_T:\cX[0,T]\to\cX_T$$
is bounded. According to Lemma 2.1, we have the following result.

\ms

\bf Corollary 2.3. \sl If {\rm(H1)$''$} holds, then
$\cA_T:\cX[0,T]\to\cX_T$ and $\cB_T:\cU[0,T]\to\cX_T$ are bounded.
Further, if {\rm(H1)$'$} holds, then
$\cA_T:L^2_\dbF(0,T;\dbR^n)\to\cX_T$ is also bounded. In particular,
all the above hold if {\rm(H1)} holds.

\rm

\ms

Recall that if $\xi\in\cX_T$, then there exists a unique $\z(\cd)\in
L^2_\dbF(0,T;\dbR^n)$ such that
$$\xi=\dbE\xi+\int_0^T\z(s)dW(s).$$
We denote
$$\z(s)=D_s\xi,\qq s\in[0,T],$$
and call it the {\it Malliavin derivative} of $\xi$
(\cite{Nualart}). Next, we have the following results which give
representation of the adjoint operators of $\cA$, $\cB$, $\cA_T$,
and $\cB_T$.

\ms

\bf Proposition 2.4. \sl The following hold:
\bel{}\left\{\ba{ll}
\ns\ds(\cA^*Y)(s)=\int_s^T\3n\(A(s)^TY(t)+\h
A(s)^T\dbE[Y(t)]+A_1(s)^TD_sY(t)+\h A_1(s)^T\dbE[D_sY(t)]\)dt,\\
\ns\ds(\cB^*Y)(s)=\int_s^T\3n\(B(s)^TY(t)+\h
B(s)^T\dbE[Y(t)]+B_1(s)^TD_sY(t)+\h B_1(s)^T\dbE[D_sY(t)]\)dt,\\
\ns\ds(\cA_T^*\xi)(s)=A(s)^T\xi+\h
A(s)^T\dbE\xi+A_1(s)^TD_s\xi+\h A_1(s)^T\dbE[D_s\xi],\\
\ns\ds(\cB_T^*\xi)(s)=B(s)^T\xi+\h B(s)^T\dbE\xi+B_1(s)^TD_s\xi+\h
B_1(s)^T\dbE[D_s\xi].\ea\right.\ee

 \ms

\it Proof. \rm For any $Y(\cd)\in L^2_\dbF(0,T;\dbR^n)$,
$$\ba{ll}
\ns\ds\lan X,\cA^*Y\ran=\lan\cA X,Y\ran=\dbE\int_0^T\lan[\cA
X](t),Y(t)\ran dt\\
\ns\ds=\1n\dbE\2n\int_0^T\3n\lan\2n\int_0^t\3n\(A(s)X(s)\1n+\1n\h
A(s)\dbE[X(s)]\)ds\2n+\2n\int_0^t\2n\(A_1(s)X(s)\1n+\1n\h
A_1(s)\dbE[X(s)]\)dW(s),Y(t)\ran dt\\
\ns\ds=\1n\dbE\2n\int_0^T\3n\int_s^T\3n\lan A(s)X(s)\1n+\1n\h
A(s)\dbE[X(s)],Y(t)\ran dtds\1n+\1n\dbE\2n\int_0^T\3n\int_s^T\3n\lan
A_1(s)X(s)\1n+\1n\h
A_1(s)\dbE[X(s)],D_sY(t)\ran dtds\\
\ns\ds=\dbE\2n\int_0^T\3n\lan X(s),\int_s^T\3n\(A(s)^TY(t)\1n+\2n\h
A(s)^T\dbE[Y(t)]\1n+\1n A_1(s)^TD_sY(t)\1n+\1n\h
A_1(s)^T\dbE[D_sY(t)]\)dt\ran ds.\ea$$
Thus, the representation of $\cA^*$ follows. Similarly, we can
obtain the representation of $\cB^*$.

\ms

Next, for any $\xi\in\cX_T$,
$$\ba{ll}
\ns\ds\lan
X,\cA_T^*\xi\ran=\lan\cA_TX,\xi\ran\\
\ns\ds=\dbE\lan\int_0^T\(A(s)X(s)+\h
A(s)\dbE[X(s)]\)ds,\xi\ran+\dbE\lan\int_0^T\(A_1(s)X(s)+\h
A_1(s)\dbE[X(s)]\)dW(s),\xi\ran\\
\ns\ds=\dbE\int_0^T\lan X(s),A(s)^T\xi+\h
A(s)^T\dbE\xi+A_1(s)^TD_s\xi+\h A_1(s)^T\dbE[D_s\xi]\ran ds.\ea$$
Therefore, the representation of $\cA_T^*$ follows. Similarly, we
can obtain the representation of $\cB_T^*$. \endpf

\ms

For completeness, let us also prove the following result.

\ms

\bf Proposition 2.5. \rm It holds
\bel{}I_T^*\xi=\xi\d_{\{T\}},\qq\forall\xi\in\cX_T,\ee
where $\d_{\{T\}}$ is the Dirac measure at $T$, and
\bel{}\dbE^*x=x^T\dbE,\qq\forall x\in\dbR^n.\ee

\ms

\it Proof. \rm First of all, since $I_T:\cX[0,T]\to\cX_T$ is
bounded, we have $I_T^*:\cX_T^*\equiv\cX_T\to\cX[0,T]^*$. For any
$\xi\in\cX_T$, and any $Y(\cd)\in\cX[0,T]$, we have
$$\lan
I_T^*\xi,Y(\cd)\ran=\lan\xi,I_TY(\cd)\ran=\dbE\lan\xi,Y(T)\ran=\dbE\int_0^T\lan
Y(s),\xi\ran\d_{T}(ds).$$
Next, since $\dbE:\cX_T\to\dbR$ is bounded, we have
$\dbE^*:\dbR\to\cX_T$. For any $\xi\in\cX_T$ and $x\in\dbR$,
$$\lan\dbE^*x,\xi\ran=\lan x,\dbE\xi\ran=\dbE\lan x,\xi\ran=x^T\dbE\xi.$$
This completes the proof. \endpf

\ms

With operators $\cA$ and $\cB$, we can write the state equation
(\ref{1.1}) as follows:
\bel{2.3}X=x+\cA X+\cB u.\ee
We now have the following result for the well-posedness of the state
equaton.

\ms

\bf Proposition 2.6. \sl Let {\rm(H1)} hold. Then for any
$(x,u(\cd))\in\dbR^n\times\cU[0,T]$, state equation $(\ref{1.1})$
admits a unique solution $X(\cd)\equiv
X(\cd\,;x,u(\cd))\in\h\cX[0,T]$.

\ms

\it Proof. \rm For any $X(\cd)\in\cX[0,T]$ and $u(\cd)\in\cU[0,T]$,
by (\ref{estimate2}), we have
$$\dbE\[\sup_{s\in[0,t]}|(\cA X)(s)|^2\]\le
\a(t)\dbE\[\sup_{s\in[0,t]}|X(s)|^2\],$$
with $\a(t)\in(0,1)$ when $t>0$ is small. Hence, by contraction
mapping theorem, we obtain well-posedness of the state equation on
$[0,t]$. Then by a usual continuation argument, we obtain the
well-posedness of the state equation on $[0,T]$. \endpf

\ms

From (\ref{estimate2}), we see that if for some $\e>0$,
\bel{2.14}\left\{\ba{ll}
\ns\ds A(\cd),\h A(\cd)\in L^{1+\e}(0,T;\dbR^{n\times n}),\qq
A_1(\cd),\h A_1(\cd)\in L^{2+\e}(0,T;\dbR^{n\times n}),\\
\ns\ds B(\cd),\h B(\cd)\in L^2(0,T;\dbR^{n\times m}),\qq B_1(\cd),\h
B_1(\cd)\in L^\infty(0,T;\dbR^{n\times m}),\ea\right.\ee
then the result of Proposition 2.6 holds. It is ready to see that
(\ref{2.14}) is stronger than (H1)$''$ and weaker than (H1)$'$. In
what follows, for convenient, we will assume (H1). However, we keep
in mind that (H1) can actually be relaxed.

\ms

Proposition 2.6 tells us that under, say, (H1), the operator
$I-\cA:\h\cX[0,T]\to\h\cX[0,T]$ is invertible and the solution $X$
to the state equation corresponding to
$(x,u(\cd))\in\dbR^n\times\cU[0,T]$ is given by
$$X=(I-\cA)^{-1}x+(I-\cA)^{-1}\cB u.$$
Note that
$$\ba{ll}
\ns\ds I_T\[(I-\cA)^{-1}x+(I-\cA)^{-1}\cB
u\]=I_TX=X(T)=x+\cA_TX+\cB_Tu\\
\ns\ds=\[I+\cA_T(I-\cA)^{-1}\]x+\[\cA_T(I-\cA)^{-1}\cB+\cB_T\]u.\ea$$
Therefore,
$$I_T(I-\cA)^{-1}=I+\cA_T(I-\cA)^{-1},\qq
I_T(I-\cA)^{-1}\cB=\cA_T(I-\cA)^{-1}\cB+\cB_T.$$
Now, let
$$\left\{\ba{ll}
\ns\ds[\cQ X(\cd)](s)=Q(s)X(s),\q s\in[0,T],\qq\forall X(\cd)\in L^2_\dbF(0,T;\dbR^n),\\
\ns\ds[\h\cQ\f(\cd)](s)=\h Q(s)\f(s),\q
s\in[0,T],\qq\forall\f(\cd)\in L^2(0,T;\dbR^n),\\
\ns\ds[\cR u(\cd)](s)=R(s)u(s),\q s\in[0,T],\qq\forall u(\cd)\in\cU[0,T],\\
\ns\ds[\h\cR\psi(\cd)](s)=\h R(s)\psi(s),\q
s\in[0,T],\qq\forall\psi(\cd)\in L^2(0,T;\dbR^m),\\
\ns\ds\cG\xi=G\xi,\qq\forall\xi\in\cX_T,\qq\h\cG x=\h Gx,\qq\forall
x\in\dbR^n.\ea\right.$$
Then the cost functional can be written as
\bel{}\3n\ba{ll}
\ns\ds J(x;u(\cd))\1n=\1n\lan\1n\cQ
X,X\1n\ran\1n+\1n\lan\1n\h\cQ\dbE X,\dbE X\1n\ran\1n+\1n\lan\1n\cR
u,u\1n\ran\1n+\1n\lan\1n\h\cR\dbE u,\dbE u\1n\ran\1n+\1n\lan\1n\cG
X(T),X(T)\1n\ran\1n+\1n\lan\1n\h\cG\dbE X(T),\dbE
X(T)\1n\ran\\
\ns\ds=\lan\cQ[(I-\cA)^{-1}x+(I-\cA)^{-1}\cB
u],(I-\cA)^{-1}x+(I-\cA)^{-1}\cB u\ran\\
\ns\ds\;+\lan\h\cQ\dbE[(I-\cA)^{-1}x+(I-\cA)^{-1}\cB
u],\dbE[(I-\cA)^{-1}x+(I-\cA)^{-1}\cB u]\ran+\lan\cR u,u\ran+\lan\h\cR\dbE u,\dbE u\ran\\
\ns\ds\;+\1n\lan\cG\{[I\1n+\1n\cA_T(I\1n-\1n\cA)^{-1}]x\1n+\1n[\cA_T(I\1n-\1n\cA)^{-1}\cB\1n+\1n\cB_T]u\},
[I\1n+\1n\cA_T(I\1n-\1n\cA)^{-1}]x\1n+\1n[\cA_T(I\1n-\1n\cA)^{-1}\cB\1n+\1n\cB_T]u\ran\\
\ns\ds\;+\1n\lan\1n\h\cG\dbE\{[I\1n+\1n\cA_T(I\1n-\1n\cA)^{-1}]x\1n+\1n[\cA_T(I\1n-\1n\cA)^{-1}\cB\2n+\1n
\cB_T]u\},\dbE\{[I\1n+\1n\cA_T(I\1n-\1n\cA)^{-1}]x\2n+\1n[\cA_T(I\1n-\1n\cA)^{-1}\cB\2n+\1n\cB_T]u\}\1n\ran\\
\ns\ds\equiv\lan\Th_2u,u\ran+2\lan\Th_1x,u\ran+\lan\Th_0x,x\ran,\ea\ee
where
$$\ba{ll}
\ns\ds\Th_2=\cR+\dbE^*\h\cR\dbE+\cB^*(I-\cA^*)^{-1}\cQ(I-\cA)^{-1}\cB+\cB^*(I-\cA^*)^{-1}\dbE^*\h\cQ
\dbE(I-\cA)^{-1}\cB\\
\ns\ds\qq+\1n[\cB^*(I\1n-\1n\cA^*)^{-1}\cA_T^*\1n+\1n\cB_T^*]\cG[\cA_T(I\1n-\1n\cA)^{-1}\cB\1n+\1n\cB_T]
\1n+\1n[\cB^*(I\1n-\1n\cA^*)^{-1}\cA_T^*\1n+\1n\cB_T^*]\dbE^*\h\cG\dbE[\cA_T(I\1n-\1n\cA)^{-1}
\cB\1n+\1n\cB_T]\\
\ns\ds\q=\cR+\dbE^*\h\cR\dbE+\cB^*(I-\cA^*)^{-1}(\cQ+\dbE^*\h\cQ\dbE)(I-\cA)^{-1}\cB\\
\ns\ds\qq\q+[\cB^*(I-\cA^*)^{-1}\cA_T^*+\cB_T^*](\cG+\dbE^*\h\cG\dbE)[\cA_T(I-\cA)^{-1}\cB+\cB_T],\\
\ns\ds\Th_1=\cB^*(I-\cA^*)^{-1}\cQ(I-\cA)^{-1}+\cB^*(I-\cA^*)^{-1}\dbE^*\h\cQ\dbE(I-\cA)^{-1}\\
\ns\ds\qq\q+[\cB^*(I-\cA^*)^{-1}\cA_T^*+\cB_T^*]\cG[I+\cA_T(I-\cA)^{-1}]\\
\ns\ds\qq\q+[\cB^*(I-\cA^*)^{-1}\cA_T^*+\cB_T^*]\dbE^*\h\cG\dbE[I+\cA_T(I-\cA)^{-1}]\\
\ns\ds\q=\cB^*(I-\cA^*)^{-1}(\cQ+\dbE^*\h\cQ\dbE)(I-\cA)^{-1}
+[\cB^*(I-\cA^*)^{-1}\cA_T^*+\cB_T^*](\cG+\dbE^*\h\cG\dbE)[I+\cA_T(I-\cA)^{-1}],\\
\ns\ds\Th_0=(I-\cA^*)^{-1}\cQ(I-\cA)^{-1}+(I-\cA^*)^{-1}\dbE^*\h\cQ
\dbE(I-\cA)^{-1}\\
\ns\ds\qq\q+[I+(I-\cA^*)^{-1}\cA_T^*]\cG[I+\cA_T(I-\cA)^{-1}]
+[I+(I-\cA^*)^{-1}\cA_T^*]\dbE^*\h\cG\dbE[I+\cA_T(I-\cA)^{-1}]\\
\ns\ds\q=(I-\cA^*)^{-1}(\cQ+\dbE^*\h\cQ\dbE)(I-\cA)^{-1}+[I+(I-\cA^*)^{-1}\cA_T^*]
(\cG+\dbE^*\h\cG\dbE)[I+\cA_T(I-\cA)^{-1}].\ea$$
Consequently, for any $u(\cd),v(\cd)\in\cU[0,T]$, and $x\in\dbR^n$,
\bel{}\ba{ll}
\ns\ds J(x;v(\cd))=J(x;u(\cd)+[v(\cd)-u(\cd)])\\
\ns\ds\qq=\lan\Th_2[u+(v-u)],u+(v-u)\ran+2\lan\Th_1x,u+(v-u)\ran+\lan\Th_0x,x\ran\\
\ns\ds\qq=\lan\Th_2u,u\ran+2\lan\Th_1x,u\ran+\lan\Th_0x,x\ran+2\lan\Th_2u+\Th_1x,v-u\ran
+\lan\Th_2(v-u),v-u\ran\\
\ns\ds\qq=J(x;u(\cd))+2\lan\Th_2u+\Th_1x,v-u\ran+\lan\Th_2(v-u),v-u\ran.\ea\ee
We now present the following result whose proof is standard, making
use of the above (see \cite{Mou-Yong 2006} for details).

\ms

\bf Proposition 2.7. \sl If $u(\cd)\mapsto J(x;u(\cd))$ admits a
minimum, then
\bel{Th2ge0}\Th_2\ge0.\ee
Conversely, if in addition to $(\ref{Th2ge0})$, one has
\bel{}\Th_1x\in\Th_2\(\cU[0,T]\),\ee
then $u(\cd)\mapsto J(x;u(\cd))$ admits a minimum $\bar
u(\cd)\in\cU[0,T]$. Further, if
\bel{Th2>0}\Th_2\ge\d I,\ee
for some $\d>0$, then for any given $x\in\dbR^n$, $u(\cd)\mapsto
J(x;u(\cd))$ admits a unique minimum.

\rm

\ms

By the definition of $\Th_2$, we see that (\ref{Th2ge0}) is implied
by the following:
\bel{2.20}\cR+\dbE^*\h\cR\dbE\ge0,\q\cQ+\dbE^*\h\cQ\dbE\ge0,\q\cG+\dbE^*\h\cG\dbE\ge0,\ee
and (\ref{Th2>0}) is implied by
\bel{2.21}\cR+\dbE^*\h\cR\dbE\ge\d
I,\q\cQ+\dbE^*\h\cQ\dbE\ge0,\q\cG+\dbE^*\h\cG\dbE\ge0,\ee
for some $\d>0$. Now, we would like to present more direct
conditions under which (\ref{Th2ge0}) and (\ref{Th2>0}) hold,
respectively.

\ms

\bf Proposition 2.8. \sl Let {\rm(H1)} and {\rm(H2)$'$} hold. Then
$(\ref{Th2ge0})$ holds. Further, if {\rm(H2)$''$} holds for some
$\d>0$, then $(\ref{Th2>0})$ holds and Problem {\rm(MF)} admits a
unique solution.

\ms

\it Proof. \rm For any $\xi\in\cX_T$,
$$\ba{ll}
\ns\ds\dbE\[\lan G\xi,\xi\ran+\lan\h
G\dbE[\xi],\dbE[\xi]\ran\]=\dbE\[\lan
G(\xi-\dbE[\xi]),\xi-\dbE[\xi]\ran +\lan(G+\h
G)\dbE[\xi],\dbE[\xi]\ran\]\ge0,\\
\ns\ds\dbE\[\lan Q(s)\xi,\xi\ran+\lan\h
Q(s)\dbE[\xi],\dbE[\xi]\ran\]=\dbE\[\lan
Q(s)(\xi-\dbE[\xi]),\xi-\dbE[\xi]\ran+\lan[Q(s)+\h
Q(s)]\dbE[\xi],\dbE[\xi]\ran\]\ge0,\ea$$
and for any $\eta\in\cU_T$,
$$\dbE\[\lan R(s)\eta,\eta\ran+\lan\h
R(s)\dbE[\eta],\dbE[\eta]\ran\]=\dbE\[\lan
R(s)(\eta-\dbE[\eta]),\eta-\dbE[\eta]\ran+\lan[R(s)+\h
R(s)]\dbE[\eta],\dbE[\eta]\ran\]\ge0.$$
Thus, (\ref{Th2ge0}) holds. Next, if (H2)$''$ holds, then
$$\ba{ll}
\ns\ds\dbE\[\lan R(s)\eta,\eta\ran+\lan\h
R(s)\dbE[\eta],\dbE[\eta]\ran\]=\dbE\[\lan
R(s)(\eta-\dbE[\eta]),\eta-\dbE[\eta]\ran+\lan[R(s)+\h
R(s)]\dbE[\eta],\dbE[\eta]\ran\]\\
\ns\ds\qq\ge\d\dbE\[|\eta-\dbE[\eta]|^2+|\dbE[\eta]|^2\]=\d\dbE|\eta|^2.\ea$$
Hence, (\ref{Th2>0}) holds. \endpf

\section{Optimality Conditions.}

In this section, we first derive a necessary condition for optimal
pair of Problem (MF).

\ms

\bf Theorem 3.1. \sl Let {\rm(H1)} and {\rm(H2)} hold. Let $(\bar
X(\cd),\bar u(\cd))$ be an optimal pair of Problem {\rm(MF)}. Then
the following mean-field backward stochastic differential equation
(MF-BSDE, for short) admits a unique adapted solution
$(Y(\cd),Z(\cd))$:
\bel{MF-BSDE1}\left\{\ba{ll}
\ns\ds dY(s)=-\(A(s)^TY(s)+A_1(s)^TZ(s)+\h A(s)^T\dbE[Y(s)]+\h
A_1(s)^T\dbE[Z(s)]\\
\ns\ds\qq\qq\qq+Q(s)\bar X(s)+\h Q(s)\dbE[\bar X(s)]\)ds+Z(s)dW(s),\qq s\in[0,T],\\
\ns\ds Y(T)=G\bar X(T)+\h G\dbE[\bar X(T)],\ea\right.\ee
such that
\bel{bar u}\ba{ll}
\ns\ds R(s)\bar u(s)+B(s)^TY(s)+B_1(s)^TZ(s)+\h R(s)\dbE[\bar
u(s)]+\h B(s)^T\dbE[Y(s)]+\h
B_1(s)^T\dbE[Z(s)]=0,\\
\ns\ds\qq\qq\qq\qq\qq\qq\qq\qq\qq\qq\qq\qq\qq s\in[0,T],\q\as\ea\ee

\ms

\it Proof. \rm Let $(\bar X(\cd),\bar u(\cd))$ be an optimal pair of
Problem (MF). For any $u(\cd)\in\cU[0,T]$, let $X(\cd)$ be the state
process corresponding to the zero initial condition and the control
$u(\cd)$. Then we have
$$\ba{ll}
\ns\ds0=\dbE\Big\{\int_0^T\[\lan Q(s)\bar X(s),X(s)\ran+\lan\h
Q(s)\dbE[\bar X(s)],\dbE[X(s)]\ran+\lan R(s)\bar
u(s),u(s)\ran\\
\ns\ds\qq\qq\qq+\lan\h R(s)\dbE[\bar u(s)],\dbE[u(s)]\ran\]ds+\lan
G\bar X(T),X(T)\ran+\lan\h G\dbE[\bar
X(T)],\dbE[X(T)]\ran\Big\}\\
\ns\ds=\dbE\Big\{\int_0^T\[\lan Q(s)\bar X(s)+\h Q(s)\dbE[\bar
X(s)],X(s)\ran+\lan R(s)\bar u(s)+\h R(s)\dbE[\bar
u(s)],u(s)\ran\]ds\\
\ns\ds\qq\qq\qq+\lan G\bar X(T)+\h G\dbE[\bar
X(T)],X(T)\ran\Big\},\ea$$
On the other hand, by \cite{Buckdahn-Li-Peng 2009}, we know that
(\ref{MF-BSDE1}) admits a unique adapted solution $(Y(\cd),Z(\cd))$.
Then
$$\ba{ll}
\ns\ds\dbE\lan X(T),G\bar X(T)+\h G\dbE[\bar X(T)]\ran=\dbE\lan
X(T),Y(T)\ran\\
\ns\ds=\dbE\Big\{\int_0^T\(\lan A(s)X(s)+B(s)u(s)+\h
A(s)\dbE[X(s)]+\h B(s)\dbE[u(s)],Y(s)\ran\\
\ns\ds\qq-\lan X(s),A(s)^TY(s)+A_1(s)^TZ(s)+\h A(s)^T\dbE[Y(s)]+\h
A_1(s)^T\dbE[Z(s)]+Q(s)\bar X(s)+\h Q(s)\dbE[\bar X(s)]\ran\\
\ns\ds\qq+\lan A_1(s)X(s)+B_1(s)u(s)+\h
A_1(s)\dbE[X(s)]+\h B_1(s)\dbE[u(s)],Z(s)\ran\)ds\Big\}\\
\ns\ds=\dbE\Big\{\int_0^T\(-\lan X(s),Q(s)\bar X(s)+\h Q(s)\dbE[\bar
X(s)]\ran\\
\ns\ds\qq\qq+\lan u(s),B(s)^TY(s)+B_1(s)^TZ(s)+\h
B(s)^T\dbE[Y(s)]+\h B_1(s)^T\dbE[Z(s)]\ran\)ds\Big\}.\ea$$
Hence,
$$\ba{ll}
\ns\ds\dbE\2n\int_0^T\2n\3n\(\2n\lan u(s),B(s)^TY(s)\1n+\1n
B_1(s)^TZ(s)\1n+\1n\h B(s)^T\dbE[Y(s)]\1n+\1n\h
B_1(s)^T\dbE[Z(s)]\1n+\1n R(s)\bar u(s)\1n+\1n\h R(s)\dbE[\bar
u(s)]\ran\2n\)ds=0,\ea$$
which leads to
$$\ba{ll}
\ns\ds R(s)\bar u(s)+B(s)^TY(s)+B_1(s)^TZ(s)+\h R(s)\dbE[\bar
u(s)]+\h B(s)^T\dbE[Y(s)]+\h B_1(s)^T\dbE[Z(s)]=0.\ea$$
This completes the proof. \endpf

\ms

From the above, we end up with the following optimality system:
(with $s$ suppressed)
\bel{MF-FBSDE1}\left\{\ba{ll}
\ns\ds d\bar X=\(A\bar X+B\bar u+\h A\dbE[\bar X]+\h B\dbE[\bar
u]\)ds+\(A_1\bar X+B_1\bar u+\h A_1\dbE[\bar X]+\h
B_1\dbE[\bar u]\)dW(s),\\
\ns\ds dY=-\(A^TY+A_1^TZ+Q\bar X+\h
A^T\dbE[Y]+\h A_1^T\dbE[Z]+\h Q\dbE[\bar X]\)ds+ZdW(s),\\
\ns\ds X(0)=x,\qq Y(T)=G\bar X(T)+\h G\dbE[\bar
X(T)],\\
\ns\ds R\bar u+\h R\dbE[\bar u]+B^TY+B_1^TZ+\h B^T\dbE[Y]+\h
B_1^T\dbE[Z]=0.\ea\right.\ee
This is called a (coupled) forward-backward stochastic differential
equations of mean-field type (MF-FBSDE, for short). Note that the
coupling comes from the last relation (which is essentially the
maximum condition in the usual Pontryagin type maximum principle).
The 4-tuple $(\bar X(\cd),\bar u(\cd),Y(\cd),Z(\cd))$ of
$\dbF$-adapted processes satisfying the above is called an {\it
adapted solution} of (\ref{MF-FBSDE1}). We now look at the
sufficiency of the above result.

\ms

\bf Theorem 3.2. \sl Let {\rm(H1)}, {\rm(H2)}, and $(\ref{Th2ge0})$
hold. Suppose $(\bar X(\cd),\bar u(\cd),Y(\cd),Z(\cd))$ is an
adapted solution to the MF-FBSDE $(\ref{MF-FBSDE1})$. Then $(\bar
X(\cd),\bar u(\cd))$ is an optimal pair.

\ms

\it Proof. \rm Let $(\bar X(\cd),\bar u(\cd),Y(\cd),Z(\cd))$ be an
adapted solution to the MF-FBSDE. For any $u(\cd)\in\cU[0,T]$, let
$$X_1(\cd)\equiv X(\cd\,;0,u(\cd)-\bar u(\cd)).$$
Then
$$X(s;x,u(\cd))=\bar X(s)+X_1(s),\qq s\in[0,T].$$
Hence, (suppressing $s$)
$$\ba{ll}
\ns\ds J(x;u(\cd))-J(x;\bar u(\cd))\\
\ns\ds=2\dbE\Big\{\int_0^T\[\lan Q\bar X,X_1\ran+\lan\h
Q\dbE[\bar X],\dbE[X_1]\ran+\lan R\bar u,u-\bar u\ran+\lan\h R\dbE[\bar u],\dbE[u-\bar u]\ran\]ds\\
\ns\ds\qq+\lan G\bar X(T),X_1(T)\ran+\lan\h G\dbE[\bar X(T)],\dbE[X_1(T)]\ran\Big\}\\
\ns\ds\qq+\dbE\Big\{\int_0^T\[\lan QX_1,X_1\ran+\lan\h
Q\dbE[X_1],\dbE[X_1]\ran+\lan R(u-\bar u),u-\bar u\ran+\lan\h R\dbE[u-\bar u],\dbE[u-\bar u]\ran\]ds\\
\ns\ds\qq+\lan GX_1(T),X_1(T)\ran+\lan\h G\dbE[X_1(T)],\dbE[X_1(T)]\ran\Big\}\\
\ns\ds=2\dbE\Big\{\int_0^T\[\lan X_1,Q\bar X+\h Q\dbE[\bar
X]\ran+\lan u-\bar u,R\bar u+\h R\dbE[\bar u]\ran\]ds+\lan
X_1(T),G\bar X(T)+\h
G\dbE[\bar X(T)]\ran\Big\}\\
\ns\ds\qq+J(0;u(\cd)-\bar u(\cd)).\ea$$
Note that
$$\ba{ll}
\ns\ds\dbE\lan X_1(T),G\bar X(T)+\h G\dbE[\bar X(T)]\ran=\dbE\lan
X_1(T),Y(T)\ran\\
\ns\ds=\dbE\Big\{\int_0^T\(\lan AX_1+B(u-\bar u)+\h
A\dbE[X_1]+\h B\dbE[u-\bar u],Y\ran\\
\ns\ds\qq-\lan X_1,A^TY+A_1^TZ+\h A^T\dbE[Y]+\h
A_1^T\dbE[Z]+Q\bar X+\h Q\dbE[\bar X]\ran\\
\ns\ds\qq+\lan A_1X_1+B_1(u-\bar u)+\h
A_1\dbE[X_1]+\h B_1\dbE[u-\bar u],Z\ran\)ds\Big\}\\
\ns\ds=\dbE\Big\{\int_0^T\(-\lan X_1,Q\bar X+\h Q\dbE[\bar
X]\ran+\lan u-\bar u,B^TY+B_1^TZ+\h B^T\dbE[Y]+\h
B_1^T\dbE[Z]\ran\)ds\Big\}.\ea$$
Thus,
$$\ba{ll}
\ns\ds J(x;u(\cd))-J(x;\bar u(\cd))\\
\ns\ds=2\dbE\int_0^T\lan u-\bar u,R\bar u+\h R\dbE[\bar
u]+B^TY+B_1^TZ+\h B^T\dbE[Y]+\h B_1^T\dbE[Z]\ran
ds+J(0;u(\cd)-\bar u(\cd))\\
\ns\ds=J(0;u(\cd)-\bar u(\cd))=\lan\Th_2(u-\bar u),u-\bar
u\ran\ge0.\ea$$
Hence, $(\bar X(\cd),\bar u(\cd))$ is optimal. \endpf

\ms

We have the following corollary.

\ms

\bf Corollary 3.3. \sl Let {\rm(H1)} and {\rm(H2)$''$} hold. Then
MF-FBSDE $(\ref{MF-FBSDE1})$ admits a unique adapted solution $(\bar
X(\cd),\bar u(\cd),Y(\cd),Z(\cd))$ of which $(\bar X(\cd),\bar
u(\cd))$ is the unique optimal pair of Problem {\rm(MF)}.

\ms

\it Proof. \rm We know from Proposition 2.8 that under (H1) and
(H2)$''$, Problem (MF) admits a unique optimal pair $(\bar
X(\cd),\bar u(\cd))$. Then by Theorem 3.1, for some
$(Y(\cd),Z(\cd))$, the 4-tuple $(\bar X(\cd),\bar
u(\cd),Y(\cd),Z(\cd))$ is an adapted solution to MF-FBSDE
(\ref{MF-FBSDE1}). Next, if (\ref{MF-FBSDE1}) has another adapted
solution $(\wt X(\cd),\wt u(\cd),\wt Y(\cd),\wt Z(\cd))$. Then by
Theorem 3.2, $(\wt X(\cd),\wt u(\cd))$ must be an optimal pair of
Problem (MF). Hence, by the uniqueness of the optimal pair of
Problem (MF), we must have
$$\wt X(\cd)=\bar X(\cd),\qq\wt u(\cd)=\bar u(\cd).$$
Then by the uniqueness of MF-BSDE (\ref{MF-BSDE1}), one must have
$$\wt Y(\cd)=Y(\cd),\qq\wt Z(\cd)=Z(\cd),$$
proving the corollary. \endpf

\ms

\section{Decoupling the MF-FBSDE and Riccati Equations}

From Corollary 3.3, we see that under (H1) and (H2)$''$, to solve
Problem (MF), we need only to solve MF-FBSDE (\ref{MF-FBSDE1}). To
solve (\ref{MF-FBSDE1}), we use the idea of decoupling inspired by
the Four-Step Scheme (\cite{MPY 1994,MY 1999}). More precisely, we
have the following result. For simplicity, we have suppressed the
time variable $s$ below.

\ms

\bf Theorem 4.1. \sl Let {\rm(H1)} and {\rm(H2)$''$} hold. Then the
unique adapted solution $(\bar X(\cd),\bar u(\cd),Y(\cd),Z(\cd))$ of
MF-FBSDE $(\ref{MF-FBSDE1})$ admits the following representation:
\bel{representation}\left\{\ba{ll}
\ns\ds\bar u=-K_0^{-1}(B^TP+B_1^TPA_1)\(\bar X-\dbE[\bar
X]\)-K_1^{-1}\[(B+\h B)^T\Pi+(B_1+\h B_1)^TP(A_1+\h A_1)\]\dbE[\bar
X],\\
\ns\ds Y=P\(\bar X-\dbE[\bar X]\)+\Pi\dbE[\bar X],\\
\ns\ds Z=\[PA_1-PB_1K_0^{-1}(B^TP+B_1^TPA_1)\]\(\bar X-\dbE[\bar
X]\)\\
\ns\ds\qq+\[P(A_1+\h A_1)-P(B_1+\h B_1)K_1^{-1}\((B+\h
B)^T\Pi+(B_1+\h B_1)^TP(A_1+\h A_1)\)\]\dbE[\bar X].\ea\right.\ee
where
\bel{K}K_0=R+B_1^TPB_1,\qq K_1=R+\h R+(B_1+\h B_1)^TP(B_1+\h
B_1).\ee
and $P(\cd)$ and $\Pi(\cd)$ are solutions to the following Riccati
equations, respectively:
\bel{Riccati P}\left\{\2n\ba{ll}
\ns\ds P'+PA+A^TP+A_1^TPA_1+Q-(PB\1n+\1n A_1^TPB_1)K_0^{-1}
(B^TP+B_1^TPA_1)=0,\q s\in[0,T],\\
\ns\ds P(T)=G,\ea\right.\ee
and
\bel{Riccati Pi}\left\{\3n\ba{ll}
\Pi'+\Pi\[(A+\h A)-(B+\h B)K_1^{-1}(B_1+\h B_1)^TP(A_1+\h A_1)\]\\
\ns\ds\qq+\[(A+\h A)^T-(A_1+\h A_1)^TP(B_1+\h B_1)K_1^{-1}(B+\h
B)^T\]\Pi-\Pi(B+\h B)K_1^{-1}
(B+\h B)^T\Pi\\
\ns\ds\qq+(A_1+\h A_1)^T\[P-P(B_1+\h B_1)K_1^{-1}(B_1+\h
B_1)^TP\](A_1+\h A_1)+Q+\h Q=0,\qq s\in[0,T],\\
\ns\ds\Pi(T)=G+\h G.\ea\right.\ee
Finally, $\bar X(\cd)$ solves the following closed-loop system
\bel{closed-loop}\left\{\ba{ll}
\ns\ds d\bar X=\Big\{\[A-BK_0^{-1}(B^TP+B_1^TPA_1)\]\(\bar X-\dbE[\bar X]\)\\
\ns\ds\qq+\[(A+\h A)-(B+\h B)K_1^{-1}\((B+\h B)^T\Pi+(B_1+\h
B_1)^TP(A_1+\h A_1)\)\]\dbE[\bar X]\Big\}dt\\
\ns\ds\qq+\Big\{\[A_1-B_1K_0^{-1}(B^TP+B_1^TPA_1)\]\(\bar X-\dbE[\bar X]\)\\
\ns\ds\qq+\[(A_1+\h A_1)-(B_1+\h B_1)K_1^{-1}\((B+\h B)^T\Pi+(B_1+\h
B_1)^TP(A_1+\h A_1)\)\]\dbE[\bar X]\Big\}dW(t),\\
\ns\ds\bar X(0)=x.\ea\right.\ee

\ms

\rm

\it Proof. \rm First of all, under (H1) and (H2)$''$, Riccati
equation (\ref{Riccati P}) admits a unique solution $P(\cd)$ which
is positive definite matrix valued (\cite{Yong-Zhou 1999}). With
such a function $P(\cd)$, $K_1$ defined in (\ref{K}) is positive
definite. Next, we note that
$$\ba{ll}
\ns\ds P-P(B_1+\h B_1)K_1^{-1}(B_1+\h B_1)^TP\\
\ns\ds\q=P\1n-\1n P(B_1\1n+\1n\h B_1)\[R\1n+\1n\h
R\1n+\1n(B_1\1n+\1n\h B_1)^TP(B_1\1n+\1n\h
B_1)\]^{-1}(B_1\1n+\1n\h B_1)^TP\equiv P\1n-\1n P\wt B(\wt R\1n+\1n\wt B^TP\wt B)^{-1}\wt B^TP\\
\ns\ds\q=P^{1\over2}\[I\1n-\1n P^{1\over2}\wt B\wt
R^{-{1\over2}}\(I\1n+\1n\wt R^{-{1\over2}}\wt
B^TP^{1\over2}P^{1\over2}\wt B\wt R^{-{1\over2}}\)^{-1}\wt
R^{-{1\over2}}\wt B^TP^{1\over2}\]P^{1\over2}\equiv
P^{1\over2}\[I\1n-\1n\G(I\1n+\1n\G^T\G)^{-1}\G^T\]P^{1\over2}\\
\ns\ds=\1n P^{1\over2}\1n(I\1n+\1n\G\G^T\1n)^{-1}\1n
P^{1\over2}\1n\equiv\1n P^{1\over2}(I\1n+\1n P^{1\over2}\wt B\wt
R^{-1}\wt B^T\1n P^{1\over2})^{-1}\1n P^{1\over2}\1n\equiv\1n
P^{1\over2}\1n[I\2n+\2n P^{1\over2}\1n(B_1\2n+\1n\h B_1)(R\1n+\1n\h
R)^{-1}\1n(B_1\2n+\1n\h B_1)^T\1n P^{1\over2}]^{-1}\1n
P^{1\over2}\1n\ge\1n0.\ea$$
In the above, we have denoted
$$\wt B=B_1+\h B_1,\q\wt R=R+\h R,\q\G=P^{1\over2}\wt B\wt R^{1\over2},$$
and used the fact
$$I-\G(I+\G^T\G)^{-1}\G^T=(I+\G^T\G)^{-1}.$$
Hence, Riccati equation (\ref{Riccati Pi}) is equivalent to the
following:
\bel{}\left\{\ba{ll}
\ns\ds\Pi'+\Pi\[(A+\h A)-(B+\h B)K_1^{-1}(B_1+\h B_1)^TP(A_1+\h A_1)\]\\
\ns\ds~+\[(A\1n+\1n\h A)^T\1n-\1n(A_1\1n+\1n\h A_1)^TP(B_1\1n+\1n\h
B_1)K_1^{-1}(B\1n+\1n\h B)^T\]\Pi\1n-\1n\Pi(B\1n+\1n\h B)K_1^{-1}
(B\1n+\1n\h B)^T\Pi\\
\ns\ds~+\1n(A_1\2n+\1n\h A_1)^T\1n P^{1\over2}\[I\1n+\1n
P^{1\over2}(B_1\2n+\1n\h B_1)(R\1n+\1n\h R)^{-1}(B_1\2n+\1n\h
B_1)^T\1n P^{1\over2}\]^{-1}\2n P^{1\over2}(A_1\2n+\1n\h A_1)\1n+\1n
Q+\1n\h Q=0,\\
\ns\ds\Pi(T)=G+\h G.\ea\right.\ee
Since
$$\ba{ll}
\ns\ds G+\h G\ge0,\\
\ns\ds(A_1\1n+\1n\h A_1)^TP^{1\over2}\[I\1n+\1n
P^{1\over2}(B_1\1n+\1n\h B_1)(R\1n+\1n\h R)^{-1}(B_1\1n+\1n\h
B_1)^TP^{1\over2}\]^{-1}P^{1\over2}(A_1\1n+\1n\h A_1)\1n+\1n Q+\1n\h
Q\ge0,\ea$$
according to \cite{Yong-Zhou 1999}, Riccati equation (\ref{Riccati
Pi}) admits a unique solution $\Pi(\cd)$ which is also positive
definite matrix valued.

\ms

Now, suppose $(\bar X(\cd),\bar u(\cd),Y(\cd),Z(\cd))$ is the
adapted solution to (\ref{MF-FBSDE1}). Assume that
\bel{Y}Y(s)=P(s)X(s)+\h P(s)\dbE[X(s)],\qq s\in[0,T],\ee
for some deterministic and differentiable functions $P(\cd)$ and $\h
P(\cd)$ such that
\bel{}P(T)=G,\qq\h P(T)=\h G.\ee
For the time being, we do not assume that $P(\cd)$ is the solution
to (\ref{Riccati P}). Then (suppressing $s$)
\bel{4.9}\ba{ll}
\ns\ds\(-A^TY-A_1^TZ-Q\bar
X-\h A^T\dbE[Y]-\h A_1^T\dbE[Z]-\h Q\dbE[\bar X]\)ds+ZdW(s)\\
\ns\ds=dY=d\(P\bar X+\h P\dbE[\bar X]\)\\
\ns\ds=\[P'\bar X+P\(A\bar X+B\bar u+\h A\dbE[\bar X]+\h B\dbE[\bar
u]\)+\h P'\dbE[\bar X]+\h P\((A+\h A)\dbE[\bar
X]+(B+\h B)\dbE[\bar u]\)\]ds\\
\ns\ds\qq\qq+P\(A_1\bar X+B_1\bar u+\h A_1\dbE[\bar X]+\h B_1\dbE[\bar u]\)dW(s)\\
\ns\ds=\[(P'+PA)\bar X+PB\bar u+\(\h P'+\h P(A+\h A)+P\h
A\)\dbE[\bar X]+\(\h P(B+\h B)+P\h B\)\dbE[\bar u]\]ds\\
\ns\ds\qq\qq+P\(A_1\bar X+B_1\bar u+\h A_1\dbE[\bar X]+\h
B_1\dbE[\bar u]\)dW(s)\ea\ee
Comparing the diffusion terms,  we should have
\bel{Z}Z=P\(A_1\bar X+B_1\bar u+\h A_1\dbE[\bar X]+\h B_1\dbE[\bar
u]\).\ee
This yields from (\ref{bar u}) that
$$\ba{ll}
\ns\ds0=R\bar u+\h R\dbE[\bar u]+B^TY+\h B^T\dbE[Y]+B_1^TZ+\h B_1^T\dbE[Z]\\
\ns\ds\q=R\bar u+\h R\dbE[\bar u]+B^T(P\bar X+\h P\dbE[\bar X])+\h
B^T(P+\h P)\dbE[\bar X]\\
\ns\ds\qq+B_1^TP\(A_1\bar X+B_1\bar u+\h A_1\dbE[\bar X]+\h
B_1\dbE[\bar u]\)+\h B_1^TP\((A_1+\h A_1)\dbE[\bar X]+(B_1
+\h B_1)\dbE[\bar u]\)\\
\ns\ds\q=(R+B_1^TPB_1)\bar u+(\h R+B_1^TP\h B_1+\h B_1^TPB_1+\h
B_1^TP\h B_1)\dbE[\bar u]+(B^TP+B_1^TPA_1)\bar X\\
\ns\ds\qq+\(B^T\h P+\h B^T(P+\h P)+B_1^TP\h A_1+\h B_1^TP(A_1+\h
A_1)\)\dbE[\bar X].\ea$$
Taking $\dbE$, we obtain
$$\ba{ll}
\ns\ds0=\[R\1n+\1n\h R\1n+\1n(B_1\1n+\1n\h B_1)^TP(B_1\1n+\1n\h
B_1)\]\dbE[\bar u]\1n+\1n\[(B\1n+\1n\h B)^T(P\1n+\1n\h
P)\1n+\1n(B_1\1n+\1n\h B_1)^TP(A_1\1n+\1n\h A_1)\]\dbE[\bar X]\\
\ns\ds\q\equiv K_1\dbE[\bar u]\1n+\1n\[(B\1n+\1n\h B)^T(P\1n+\1n\h
P)\1n+\1n(B_1\1n+\1n\h B_1)^TP(A_1\1n+\1n\h A_1)\]\dbE[\bar X].\ea$$
Assuming
$$K_1\equiv K_1(P)\deq R+\h R+(B_1+\h B_1)^TP(B_1+\h B_1)$$
to be invertible, one gets
$$\dbE[\bar u]=-K_1^{-1}\[(B+\h B)^T(P+\h P)+(B_1+\h B_1)^TP(A_1+\h
A_1)\]\dbE[\bar X].$$
Then we have
$$\ba{ll}
\ns\ds0=(R+B_1^TPB_1)\bar u+\(\h R+B_1^TP\h B_1+\h B_1^TPB_1+\h
B_1^TP\h B_1\)\dbE[\bar u]+(B^TP+B_1^TPA_1)\bar X\\
\ns\ds\qq+\(B^T\h P+\h B^T(P+\h P)+B_1^TP\h A_1+\h B_1^TP(A_1+\h
A_1)\)\dbE[\bar X]\\
\ns\ds\q\equiv K_0\bar u+(B^TP+B_1^TPA_1)\bar X++\(B^T\h P+\h
B^T(P+\h P)+B_1^TP\h A_1+\h B_1^TP(A_1+\h
A_1)\)\dbE[\bar X]\\
\ns\ds\qq-\(\h R\1n+\1n B_1^TP\h B_1\1n+\1n\h B_1^TPB_1\1n+\1n\h
B_1^TP\h B_1\)K_1^{-1}\[(B\1n+\1n\h B)^T(P\1n+\1n\h
P)\1n+\1n(B_1\1n+\1n\h B_1)^TP(A_1\1n+\1n\h A_1)\]\dbE[\bar X].\ea$$
Consequently, by assuming
$$K_0\equiv K_0(P)\deq R+B_1^TPB_1$$
to be invertible, we obtain
\bel{u}\ba{ll}
\ns\ds\bar u=-K_0^{-1}(B^TP+B_1^TPA_1)\bar X-K_0^{-1}\(B^T\h P+\h
B^T(P+\h P)+B_1^TP\h A_1+\h
B_1^TP(A_1+\h A_1)\)\dbE[\bar X]\\
\ns\ds\qq+K_0^{-1}\(\1n\h R\1n+\1n B_1^T\1n P\h B_1\1n+\1n\h
B_1^TPB_1\1n+\1n\h B_1^T\1n P\h B_1\1n\)K_1^{-1}\[(B\1n+\1n\h
B)^T\1n(P\1n+\1n\h P)\1n+\1n(B_1\1n+\1n\h B_1)^T\1n P(A_1\1n+\1n\h
A_1)\]\dbE[\bar X]\\
\ns\ds=-K_0^{-1}(B^TP+B_1^TPA_1)\bar X-K_0^{-1}\(B^T\h P+\h B^T(P+\h
P)+B_1^TP\h A_1+\h B_1^TP(A_1+\h A_1)\)\dbE[\bar
X]\\
\ns\ds\qq+K_0^{-1}\(K_1-K_0\)K_1^{-1}\[(B+\h B)^T(P+\h P)+(B_1+\h
B_1)^TP(A_1+\h A_1)\]\dbE[\bar X]\\
\ns\ds=-K_0^{-1}(B^TP+B_1^TPA_1)\bar X-K_0^{-1}\(B^T\h P+\h B^T(P+\h
P)+B_1^TP\h A_1+\h B_1^TP(A_1+\h A_1)\)\dbE[\bar
X]\\
\ns\ds\qq+K_0^{-1}\[(B+\h B)^T(P+\h P)+(B_1+\h B_1)^TP(A_1+\h
A_1)\]\dbE[\bar X]\\
\ns\ds\qq-K_1^{-1}\[(B+\h B)^T(P+\h P)+(B_1+\h B_1)^TP(A_1+\h
A_1)\]\dbE[\bar X]\\
\ns\ds=-K_0^{-1}(B^TP+B_1^TPA_1)\bar X+K_0^{-1}\[B^TP+B_1^TPA_1\]\dbE[\bar X]\\
\ns\ds\qq-K_1^{-1}\[(B+\h B)^T(P+\h P)+(B_1+\h B_1)^TP(A_1+\h
A_1)\]\dbE[\bar X].\ea\ee
Here, we note that
$$\ba{ll}
\ns\ds K_1-K_0=R+\h R+(B_1+\h B_1)^TP(B_1+\h B_1)-R-B_1^TPB_1\\
\ns\ds\qq\qq=\h R+B_1^TP\h B_1+\h B_1^TPB_1+\h B_1^TP\h B_1.\ea$$
Hence, comparing the drift terms in (\ref{4.9}), we have
$$\ba{ll}
\ns\ds0=(P'+PA)\bar X+PB\bar u+\(\h P'+\h P(A+\h A)+P\h
A\)\dbE[\bar X]+\(\h P(B+\h B)+P\h B\)\dbE[\bar u]\\
\ns\ds\qq+A^T\(P\bar X+\h P\dbE[\bar X]\)+A_1^TP\(A_1\bar X+B_1\bar
u+\h A_1\dbE[\bar X]+\h B_1\dbE[\bar u]\)\\
\ns\ds\qq+\h A^T(P+\h P)\dbE[\bar X]+\h A_1^TP\[(A_1+\h
A_1)\dbE[\bar X]+(B_1+\h B_1)\dbE[\bar u]\] +Q\bar X+\h Q\dbE[\bar X]\\
%
%\ns\ds=\[P'+PA+A^TP+A_1^TPA_1+Q\]\bar X+(PB+A_1^TPB_1)\bar u\\
%
%\ns\ds\qq+\[\h P'+\h P(A+\h A)+P\h A+A^T\h P+A_1^TP\h A_1+\h
%A^T(P+\h P)+\h A_1^TP(A_1+\h A_1)+\h Q\]\dbE[\bar X]\\
%
%\ns\ds\qq+\[\h P(B+\h B)+P\h B+A_1^TP\h B_1+\h A_1^TP(B_1+\h
%B_1)\]\dbE[\bar u]\\
%
\ns\ds=\[P'+PA+A^TP+A_1^TPA_1+Q\]\bar X+(PB+A_1^TPB_1)\bar u\\
\ns\ds\qq+\[\h P'\1n+\1n\h P(A\1n+\1n\h A)\1n+\1n(A\1n+\1n\h A)^T\h
P\1n+\1n P\h A\1n+\1n\h
A^TP\1n+\1n(A_1\1n+\1n\h A_1)^TP(A_1\1n+\1n\h A_1)\1n-\1n A_1^TPA_1\1n+\1n\h Q\]\dbE[\bar X]\\
\ns\ds\qq+\[\h P(B+\h B)+P\h B+A_1^TP\h B_1+\h A_1^TP(B_1+\h
B_1)\]\dbE[\bar u]\\
\ns\ds=\[P'+PA+A^TP+A_1^TPA_1+Q\]\bar X\\
\ns\ds\qq+(PB+A_1^TPB_1)\Big\{-K_0^{-1}(B^TP+B_1^TPA_1)\bar X+K_0^{-1}\[B^TP+B_1^TPA_1\]\dbE[\bar X]\\
\ns\ds\qq-K_1^{-1}\[(B+\h B)^T(P+\h P)+(B_1+\h B_1)^TP(A_1+\h
A_1)\]\dbE[\bar X]\Big\}\\
\ns\ds\qq+\[\h P'\1n+\1n\h P(A\1n+\1n\h A)\1n+\1n(A\1n+\1n\h A)^T\h
P\1n+\1n P\h A\1n+\1n\h
A^TP\1n+\1n(A_1\1n+\1n\h A_1)^TP(A_1\1n+\1n\h A_1)\1n-\1n A_1^TPA_1\1n+\1n\h Q\]\dbE[\bar X]\\
\ns\ds\qq+\[\h P(B+\h B)+P\h B+A_1^TP\h B_1+\h A_1^TP(B_1+\h
B_1)\]\\
\ns\ds\qq\qq\cd\Big\{-K_1^{-1}\[(B+\h B)^T(P+\h P)+(B_1+\h
B_1)^TP(A_1+\h
A_1)\]\dbE[\bar X]\Big\}\\
%
%\ns\ds=\[P'+PA+A^TP+A_1^TPA_1+Q-(PB+A_1^TPB_1)K_0^{-1}(B^TP+B_1^TPA_1)\]\bar X\\
%
%\ns\ds\qq+\Big\{(PB+A_1^TPB_1)K_0^{-1}(B^TP+B_1^TPA_1)\\
%
%\ns\ds\qq-(PB+A_1^TPB_1)K_1^{-1}\[(B+\h B)^T(P+\h P)+(B_1+\h B_1)^TP(A_1+\h A_1)\]\\
%
%\ns\ds\qq+\[\h P'+\h P(A+\h A)+(A+\h A)^T\h P+P\h A+\h
%A^TP+(A_1+\h A_1)^TP(A_1+\h A_1)-A_1^TPA_1+\h Q\]\\
%
%\ns\ds\qq-\[\h P(B+\h B)+P\h B+A_1^TP\h B_1+\h A_1^TP(B_1+\h B_1)\]\\
%
%\ns\ds\qq\qq\cd K_1^{-1}\[(B+\h B)^T(P+\h P)+(B_1+\h B_1)^TP(A_1+\h A_1)\]\Big\}\dbE[\bar X]\\
%
\ns\ds=\[P'+PA+A^TP+A_1^TPA_1+Q-(PB+A_1^TPB_1)K_0^{-1}(B^TP+B_1^TPA_1)\]\bar X\\
\ns\ds\qq+\Big\{\h P'+\h P(A+\h A)+(A+\h A)^T\h P+P\h A+\h
A^TP+(A_1+\h A_1)^TP(A_1+\h A_1)-A_1^TPA_1+\h Q\\
\ns\ds\qq+(PB+A_1^TPB_1)K_0^{-1}(B^TP+B_1^TPA_1)\\
\ns\ds\qq-\[(P\1n+\1n\h P)(B\1n+\1n\h B)\1n+\1n(A_1\1n+\1n\h
A_1)^TP(B_1\1n+\1n\h B_1)\]K_1^{-1}\[(B\1n+\1n\h B)^T(P\1n+\1n\h
P)\1n+\1n(B_1\1n+\1n\h B_1)^TP(A_1\1n+\1n\h A_1)\]\1n\Big\}\dbE[\bar
X].\ea$$
Therefore, by choosing $P(\cd)$ to be the solution to Riccati
equation (\ref{Riccati P}), we have that $K_0$ and $K_1$ are
positive definite, and the above leads to
$$\ba{ll}
\ns\ds\Big\{\h P'+\h P(A+\h A)+(A+\h A)^T\h P+P\h A+\h
A^TP+(A_1+\h A_1)^TP(A_1+\h A_1)-A_1^TPA_1+\h Q\\
\ns\ds+(PB+A_1^TPB_1)K_0^{-1}(B^TP+B_1^TPA_1)\\
\ns\ds-\[(P\1n+\1n\h P)(B\1n+\1n\h B)\1n+\1n(A_1\1n+\1n\h
A_1)^TP(B_1\1n+\1n\h B_1)\]K_1^{-1}\[(B\1n+\1n\h B)^T(P\1n+\1n\h
P)\1n+\1n(B_1\1n+\1n\h B_1)^TP(A_1\1n+\1n\h A_1)\]\1n\Big\}\dbE[\bar
X]=0.\ea$$
Now, if $\h P(\cd)$ satisfies the following:
$$\ba{ll}
\ns\ds0=\h P'+\h P(A+\h A)+(A+\h A)^T\h P+P\h A+\h
A^TP+(A_1+\h A_1)^TP(A_1+\h A_1)-A_1^TPA_1+\h Q\\
\ns\ds\qq+(PB+A_1^TPB_1)K_0^{-1}(B^TP+B_1^TPA_1)\\
\ns\ds\qq-\[(P+\h P)(B+\h B)+(A_1+\h A_1)^TP(B_1+\h
B_1)\]K_1^{-1}\[(B+\h B)^T(P+\h P)+(B_1+\h B_1)^TP(A_1+\h
A_1)\]\\
\ns\ds\q=\h P'+\h P\[A+\h A-(B+\h B)K_1^{-1}\((B+\h B)^TP+(B_1+\h B_1)^TP(A_1+\h A_1)\)\]\\
\ns\ds\qq+\[(A+\h A)^T-\(P(B+\h B)+(A_1+\h A_1)^TP(B_1+\h B_1)\)K_1^{-1}(B+\h B)^T\]\h P
-\h P(B+\h B)K_1^{-1}(B+\h B)^T\h P\\
\ns\ds\qq+P\h A+\h
A^TP+(A_1+\h A_1)^TP(A_1+\h A_1)-A_1^TPA_1+\h Q+(PB+A_1^TPB_1)K_0^{-1}(B^TP+B_1^TPA_1)\\
\ns\ds\qq-\[P(B+\h B)+(A_1+\h A_1)^TP(B_1+\h B_1)\]K_1^{-1}\[(B+\h
B)^TP+(B_1+\h B_1)^TP(A_1+\h A_1)\],\ea$$
then $(Y(\cd),Z(\cd))$ defined by (\ref{Y}) and (\ref{Z}) with $\bar
u(\cd)$ given by (\ref{u}) satisfies the MF-BSDE in
(\ref{MF-FBSDE1}). Hence, we introduce the following Riccati
equation for $\h P(\cd)$:
\bel{Riccati hP}\left\{\ba{ll}
\ns\ds\h P'+\h P\[A+\h A-(B+\h B)K_1^{-1}\((B+\h B)^TP+(B_1+\h B_1)^TP(A_1+\h A_1)\)\]\\
\ns\ds\q+\[(A+\h A)^T-\(P(B+\h B)+(A_1+\h A_1)^TP(B_1+\h
B_1)\)K_1^{-1}(B+\h B)^T\]\h P\\
\ns\ds\q -\h P(B+\h B)K_1^{-1}(B+\h B)^T\h P+P\h A+\h
A^TP+(A_1+\h A_1)^TP(A_1+\h A_1)-A_1^TPA_1+\h Q\\
\ns\ds\q+(PB+A_1^TPB_1)K_0^{-1}(B^TP+B_1^TPA_1)\\
\ns\ds\q-\[P(B+\h B)+(A_1+\h A_1)^TP(B_1+\h B_1)\]K_1^{-1}\[(B+\h
B)^TP+(B_1+\h B_1)^TP(A_1+\h A_1)\]=0,\\
\ns\ds\h P(T)=\h G.\ea\right.\ee
The solvability of this Riccati equation is not obvious since $\h G$
is just assumed to be symmetric, and
$$\ba{ll}
\ns\ds\wt Q\equiv P\h A+\h
A^TP+(A_1+\h A_1)^TP(A_1+\h A_1)-A_1^TPA_1+\h Q+(PB+A_1^TPB_1)K_0^{-1}(B^TP+B_1^TPA_1)\\
\ns\ds\qq-\[P(B+\h B)+(A_1+\h A_1)^TP(B_1+\h B_1)\]K_1^{-1}\[(B+\h
B)^TP+(B_1+\h B_1)^TP(A_1+\h A_1)\]\ea$$
is also just symmetric. To look at the solvability of such a Riccati
equation, we let
$$\Pi=P+\h P.$$
Then
$$\Pi(T)=G+\h G\ge0,$$
and consider the following
$$\ba{ll}
\ns\ds0=P'+PA+A^TP+A_1^TPA_1+Q-(PB+A_1^TPB_1)K_0^{-1}(B^TP+B_1^TPA_1)\\
\ns\ds\qq+\h P'+\h P(A+\h A)+(A+\h A)^T\h P+P\h A+\h
A^TP+(A_1+\h A_1)^TP(A_1+\h A_1)-A_1^TPA_1+\h Q\\
\ns\ds\qq+(PB+A_1^TPB_1)K_0^{-1}(B^TP+B_1^TPA_1)\\
\ns\ds\qq-\[(P\1n+\1n\h P)(B\1n+\1n\h B)\1n+\1n(A_1\1n+\1n\h
A_1)^TP(B_1\1n+\1n\h B_1)\]K_1^{-1}\[(B\1n+\1n\h B)^T(P\1n+\1n\h
P)\1n+\1n(B_1\1n+\1n\h B_1)^TP(A_1\1n+\1n\h
A_1)\]\\
\ns\ds\q=\Pi'\1n+\1n P(A\1n+\1n\h A)\1n+\1n(A\1n+\1n\h A)^TP\1n+\1n Q\1n+\1n\h Q\1n+\1n(\Pi\1n-\1n P)(A\1n+\1n
\h A)\1n+\1n(A\1n+\1n\h A)^T(\Pi\1n-\1n P)\1n+\1n(A_1\1n+\1n\h A_1)^TP(A_1\1n+\1n\h A_1)\\
\ns\ds\qq-\[\Pi(B+\h B)+(A_1+\h A_1)^TP(B_1+\h B_1)\]K_1^{-1}\[(B+\h
B)^T\Pi+(B_1+\h B_1)^TP(A_1+\h A_1)\]\\
\ns\ds\q=\Pi'+\Pi\[(A+\h A)-(B+\h B)K_1^{-1}(B_1+\h B_1)^TP(A_1+\h A_1)\]\\
\ns\ds\qq+\[(A+\h A)^T-(A_1+\h A_1)^TP(B_1+\h B_1)K_1^{-1}(B+\h
B)^T\]\Pi-\Pi(B+\h B)K_1^{-1}
(B+\h B)^T\Pi\\
\ns\ds\qq+(A_1+\h A_1)^T\[P-P(B_1+\h B_1)K_1^{-1}(B_1+\h
B_1)^TP\](A_1+\h A_1)+Q+\h Q.\ea$$
Thus, $\Pi(\cd)$ is the solution to Riccati equation (\ref{Riccati
Pi}). Consequently, Riccati equation (\ref{Riccati hP}) admits a
unique solution $\h P(\cd)=\Pi(\cd)-P(\cd)$. Then we obtain (from
(\ref{Y}) and (\ref{u}))
$$Y=P\(\bar X-\dbE[\bar X]\)+\Pi\dbE[\bar X],$$
and
$$\bar u=-K_0^{-1}(B^TP+B_1^TPA_1)\(\bar X-\dbE[\bar
X]\)-K_1^{-1}\[(B+\h B)^T\Pi+(B_1+\h B_1)^TP(A_1+\h A_1)\]\dbE[\bar
X].$$
Also, from (\ref{Z}), it follows that
$$\3n\ba{ll}
\ns\ds Z=P\(A_1\bar X+B_1\bar u+\h A_1\dbE[\bar X]+\h B_1\dbE[\bar
u]\)\\
\ns\ds\q=\1n PA_1\bar X\1n-\1n PB_1\Big\{\1n K_0^{-1}(B^T\1n
P\1n+\1n B_1^T\1n PA_1)\(\bar X\1n-\1n\dbE[\bar X]\)\2n+\1n
K_1^{-1}\[(B\2n+\1n\h B)^T\Pi\1n+\1n(B_1\1n+\1n\h B_1)^T\1n
P(A_1\2n+\1n\h
A_1)\]\dbE[\bar X]\1n\Big\}\\
\ns\ds\qq+P\h A_1\dbE[\bar X]\1n-\1nP\h B_1K_1^{-1}\[(B\1n+\1n\h
B)^T\Pi\1n+\1n(B_1\1n+\1n\h
B_1)^TP(A_1\1n+\1n\h A_1)\]\dbE[\bar X]\\
\ns\ds\q=\[PA_1-PB_1K_0^{-1}(B^TP+B_1^TPA_1)\]\(\bar X-\dbE[\bar
X]\)\\
\ns\ds\qq+\[P(A_1+\h A_1)-P(B_1+\h B_1)K_1^{-1}\((B+\h
B)^T\Pi+(B_1+\h B_1)^TP(A_1+\h A_1)\)\]\dbE[\bar X].\ea$$
Hence, (\ref{representation}) follows. Plugging the above
representations into the state equation, we obtain
$$\ba{ll}
\ns\ds d\bar X\1n=\1n\Big\{A\bar X\1n-\1n B\[K_0^{-1}(B^T\1n
P\1n+\1n B_1^TPA_1)\(\bar X\1n-\1n\dbE[\bar X]\)\1n-\1n
K_1^{-1}\((B\1n+\1n\h B)^T\Pi\1n+\1n(B_1\1n+\1n\h
B_1)^TP(A_1\1n+\1n\h A_1)\)\dbE[\bar X]\]\\
\ns\ds\qq+\h A\dbE[\bar X]-\h BK_1^{-1}\[(B+\h B)^T\Pi+(B_1+\h
B_1)^TP(A_1+\h A_1)\]\dbE[\bar X]\Big\}dt\\
\ns\ds\qq+\Big\{A_1\bar X\1n-\1n B_1\[K_0^{-1}(B^T\1n P\1n+\1n
B_1^T\1n PA_1)\(\bar X\1n-\1n\dbE[\bar X]\)\1n-\1n
K_1^{-1}\((B\1n+\1n\h B)^T\Pi\1n+\1n(B_1\1n+\1n\h B_1)^T\1n
P(A_1\1n+\1n\h A_1)\)\dbE[\bar X]\]\\
\ns\ds\qq+\h A_1\dbE[\bar X]-\h B_1K_1^{-1}\[(B+\h B)^T\Pi+(B_1+\h
B_1)^TP(A_1+\h A_1)\]\dbE[\bar X]\Big\}dW(t)\\
\ns\ds\q=\Big\{\[A-BK_0^{-1}(B^TP+B_1^TPA_1)\]\(\bar X-\dbE[\bar X]\)\\
\ns\ds\qq+\[(A+\h A)-(B+\h B)K_1^{-1}\((B+\h B)^T\Pi+(B_1+\h
B_1)^TP(A_1+\h A_1)\)\]\dbE[\bar X]\Big\}dt\\
\ns\ds\qq+\Big\{\[A_1-B_1K_0^{-1}(B^TP+B_1^TPA_1)\]\(\bar X-\dbE[\bar X]\)\\
\ns\ds\qq+\[(A_1+\h A_1)-(B_1+\h B_1)K_1^{-1}\((B+\h B)^T\Pi+(B_1+\h
B_1)^TP(A_1+\h A_1)\)\]\dbE[\bar X]\Big\}dW(t).\ea$$
This gives the closed-loop system (\ref{closed-loop}). From the
above derivation, we see that if $\bar X(\cd)$ is a solution to
(\ref{closed-loop}), and $(\bar u(\cd),Y(\cd),Z(\cd))$ are given by
(\ref{representation}), then $(\bar X(\cd),\bar
u(\cd),Y(\cd),Z(\cd))$ is an adapted solution to MF-FBSDE
(\ref{MF-FBSDE1}). By Corollary 3.3, we know that such a constructed
4-tuple $(\bar X(\cd),\bar u(\cd),Y(\cd),Z(\cd))$ is the unique
solution to (\ref{MF-FBSDE1}). \endpf

\ms

The following is a direct verification of optimality of state
feedback control.

\ms

\bf Theorem 4.2. \sl Let {\rm(H1)} and {\rm(H2)$''$} hold. Let
$P(\cd)$ and $\Pi(\cd)$ be the solutions to the Riccati equations
$(\ref{Riccati P})$ and $(\ref{Riccati Pi})$, respectively. Then the
state feedback control $\bar u(\cd)$ given in
$(\ref{representation})$ is the optimal control of Problem
{\rm(MF)}. Moreover, the optimal value of the cost is given by
\bel{}\inf_{u(\cd)\in\cU[0,T]}J(x;u(\cd))=\lan\Pi(0)x,x\ran,\qq\forall
x\in\dbR^n.\ee

\ms

\it Proof. \rm Let $P(\cd)$ and $\Pi(\cd)$ be the solutions to the
Riccati equations (\ref{Riccati P}) and (\ref{Riccati Pi}),
respectively and denote $K_0$ and $K_1$ as in (\ref{K}), which are
positive definite. Let $\h P(\cd)=\Pi(\cd)-P(\cd)$. Then $\h P(\cd)$
solves (\ref{Riccati hP}). Now, we observe
$$\ba{ll}
\ns\ds J(x;u(\cd))-\lan\Pi(0)x,x\ran=J(x;u(\cd))-\lan[P(0)+\h
P(0)]x,x\ran\\
\ns\ds=\dbE\int_0^T\Big\{\lan QX,X\ran+\lan\h
Q\dbE[X],\dbE[X]\ran+\lan Ru,u\ran+\lan
\h R\dbE[u],\dbE[u]\ran\\
\ns\ds\qq+\lan P'X,X\ran+2\lan P(AX+Bu+\h A\dbE[X]+\h
B\dbE[u]),X\ran\\
\ns\ds\qq+\lan
P(A_1X+B_1u+\h A_1\dbE[X]+\h B_1\dbE[u]),A_1X+B_1u+\h A_1\dbE[X]+\h B_1\dbE[u]\ran\\
\ns\ds\qq+\lan\h P'\dbE[X],\dbE[X]\ran+2\lan\h
P\dbE[X],(A+\h A)\dbE[X]+(B+\h B)\dbE[u]\ran\Big\}ds\\
\ns\ds=\dbE\2n\int_0^T\2n\Big\{\2n\lan(P'\1n+\1n PA\1n+\1n A^T\1n
P\1n+\1n A_1^T\1n PA_1\2n+\1n Q)X,X\ran\1n+2\lan u,(B^T\1n P\1n+\1n
B_1^T\1n PA_1)X\ran\1n+\1n\lan
(R\1n+\1n B_1^T\1n PB_1)u,u\ran\\
\ns\ds\qq+\lan\h Q\dbE[X],\dbE[X]\ran+\lan\h
R\dbE[u],\dbE[u]\ran+2\lan P(\h A\dbE[X]+\h
B\dbE[u]),\dbE[X]\ran\\
\ns\ds\qq+2\lan\1n P(A_1\dbE[X]\1n+\1n B_1\dbE[u]),\h
A_1\dbE[X]\1n+\1n\h B_1\dbE[u]\1n\ran\1n+\1n\lan\1n P(\h
A_1\dbE[X]\1n+\1n\h B_1\dbE[u]),\h A_1\dbE[X]\1n+\1n\h
B_1\dbE[u]\1n\ran\\
\ns\ds\qq+\lan\h P'\dbE[X],\dbE[X]\ran+2\lan\h
P\dbE[X],(A+\h A)\dbE[X]+(B+\h B)\dbE[u]\ran\Big\}ds\\
\ns\ds=\dbE\int_0^T\3n\Big\{\2n\lan(P'\1n+\1n PA\1n+\1n A^TP\1n+\1n
A_1^TPA_1\1n+\1n Q)X,X\ran\1n+2\lan
u,(B^TP\1n+\1n B_1^TPA_1)X\ran\1n+\1n\lan K_0u,u\ran\\
\ns\ds\qq+\lan[\h P'\1n+\1n P\h A\1n+\1n\h A^T\1n P\1n+\1n\h
A_1^T\1n PA_1\1n+\1n A_1^T\1n P\h A_1\1n+\1n\h A_1^T\1n P\h
A_1\1n+\1n\h P(A\1n+\1n\h A)\1n+\1n(A\1n+\1n\h A)^T\1n\h P\1n+\1n\h
Q]\dbE[X],\dbE[X]\ran\\
\ns\ds\qq+\lan(\h R+\h B_1^TPB_1+B_1^TP\h B_1+\h B_1^TP\h
B_1)\dbE[u],\dbE[u]\ran\\
\ns\ds\qq+2\lan\dbE[u],[\h B^TP+\h B_1^TPA_1+(B_1+\h B_1)^TP\h
A_1+(B+\h B)^T\h P]\dbE[X]\ran
\Big\}ds\\
\ns\ds=\dbE\int_0^T\3n\Big\{\2n\lan(P'\1n+\1n PA\1n+\1n A^TP\1n+\1n
A_1^TPA_1\1n+\1n Q)X,X\ran\1n+2\lan
u\1n-\1n\dbE[u],(B^TP\1n+\1n B_1^TPA_1)(X\1n-\1n\dbE[X])\ran\\
\ns\ds\qq\qq+\lan K_0(u-\dbE[u]),u-\dbE[u]\ran\\
\ns\ds\qq+\lan\(\h P'+P\h A+\h A^TP+\h A_1^TPA_1+A_1^TP\h A_1+\h
A_1^TP\h A_1\\
\ns\ds\qq\qq\qq+\h P(A+\h A)+(A+\h A)^T\h P+\h
Q\)\dbE[X],\dbE[X]\ran+\lan K_1\dbE[u],\dbE[u]\ran\\
\ns\ds\qq+2\lan\dbE[u],(B+\h B)^T(P+\h P)+(B_1+\h B_1)^TPA_1+(B_1+\h
B_1)^TP\h A_1\)\dbE[X]\ran \Big\}ds\ea$$
$$\ba{ll}
\ns\ds=\dbE\int_0^T\2n\Big\{\1n\lan[P'\1n+\1n PA\1n+\1n A^TP\1n+\1n
A_1^TPA_1\1n+\1n Q\1n-\1n
(PB\1n+\1n A_1PB_1)K_0^{-1}(B^TP\1n+\1n B_1^TPA_1)]X,X\ran\\
\ns\ds\qq+\Big|K_0^{1\over2}\[u-\dbE[u]+K_0^{-1}(B^TP+B_1^TPA_1)(X-\dbE[X])\]\Big|^2\\
\ns\ds\qq+\lan\[\h P'+P\h A+\h A^TP+\h A_1^TPA_1+A_1^TP\h A_1+\h
A_1^TP\h A_1+\h P(A+\h A)+(A+\h A)^T\h P+\h Q\\
\ns\ds\qq+(PB+A_1PB_1)K_0^{-1}(B^TP+B_1^TPA_1)\\
\ns\ds\qq-\((P+\h P)(B+\h B)+A_1^TP(B_1+\h B_1)+\h A_1^TP(B_1+\h
B_1)\)K_1^{-1}\\
\ns\ds\qq\q\cd\((B+\h B)^T(P+\h P)+(B_1+\h B_1)^TPA_1+(B_1+\h
B_1)^TP\h A_1\)\]\dbE[X],\dbE[X]\ran\\
\ns\ds\qq+\Big|K_1^{1\over2}\[\dbE[u]\1n+\1n K_1^{-1}\((B\1n+\1n\h
B)^T(P\1n+\1n\h P)\1n+\1n(B_1\1n+\1n\h
B_1)^TPA_1\1n+\1n(B_1\1n+\1n\h B_1)^T P\h
A_1\)\dbE[X]\]\Big|^2\Big\}ds\\
\ns\ds=\dbE\int_0^T\Big\{\Big|K_0^{1\over2}\[u-\dbE[u]+K_0^{-1}(B^TP+B_1^TPA_1)(X-\dbE[X])\]\Big|^2\\
\ns\ds\qq+\Big|K_1^{1\over2}\[\dbE[u]+K_1^{-1}\((B+\h
B)^T\Pi+(B_1+\h B_1)^TPA_1+(B_1+\h B_1)^TP\h
A_1\)\dbE[X]\]\Big|^2\Big\}ds\ge0.\ea$$
Then our claim follows. \endpf

\ms

We see that Riccati equation (\ref{Riccati Pi}) can be written as
follows:
\bel{Riccati Pi2}\left\{\3n\ba{ll}
\Pi'+\Pi(A+\h A)+(A+\h A)^T\Pi+(A_1+\h A_1)^TP(A_1+\h A_1)+Q+\h Q\\
\ns\ds-\[\Pi(B+\h B)+(A_1+\h A_1)^TP(B_1+\h B_1)\]K_1^{-1}
\[(B+\h B)^T\Pi+(B_1+\h B_1)^TP(A_1+\h A_1)\]=0,\\
\ns\ds\qq\qq\qq\qq\qq\qq\qq\qq\qq\qq s\in[0,T],\\
\ns\ds\Pi(T)=G+\h G.\ea\right.\ee
%
%\bel{Riccati P3}\left\{\2n\ba{ll}
%
%\ns\ds P'+PA+A^TP+A_1^TPA_1+Q-(PB\1n+\1n A_1^TPB_1)K_0^{-1}
%(B^TP+B_1^TPA_1)=0,\q s\in[0,T],\\
%
%\ns\ds P(T)=G,\ea\right.\ee
%
When
$$\h A=\h A_1=0,\qq\h B=\h B_1=0,\qq\h Q=0,\qq\h R=0,\qq\h G=0,$$
we have
$$K_0=K_1=R+B_1^TPB_1,$$
and the Riccati equation for $\Pi(\cd)$ can be written as
$$\left\{\ba{ll}
\ns\ds\Pi'+\Pi A+A^T\Pi-(\Pi
B+A_1^TPB_1)K_0^{-1}(B^T\Pi+B_1^TPA_1)+A_1^TPA_1+Q=0,\\
\ns\ds\Pi(T)=G.\ea\right.$$
Then
$$\ba{ll}
\ns\ds0=(\Pi-P)'+(\Pi-P)A+A^T(\Pi-P)-(\Pi
B+A_1^TPB_1)K_0^{-1}(B^T\Pi+B_1^TPA_1)\\
\ns\ds\qq+(PB+A_1^TPB_1)K_0^{-1}(B^TP+B_1^TPA_1)\\
\ns\ds\q=(\Pi-P)'+(\Pi-P)A+A^T(\Pi-P)-(\Pi-P)
BK_0^{-1}(B^T\Pi+B_1^TPA_1)\\
\ns\ds\qq-(PB+A_1^TPB_1)K_0^{-1}B^T(\Pi-P).\ea$$
Therefore, by uniqueness, we have
$$\Pi=P.$$
Consequently, the feedback control can be written as
$$\ba{ll}
\ns\ds\bar u=-K_0^{-1}(B^TP+B_1^TPA_1)\(\bar X-\dbE[\bar
X]\)-K_1^{-1}\[(B+\h B)^T\Pi+(B_1+\h B_1)^TP(A_1+\h A_1)\]\dbE[\bar
X],\\
\ns\ds\q=-K_0^{-1}(B^TP+B_1^TPA_1)\(\bar X-\dbE[\bar
X]\)-K_0^{-1}\(B^TP+B_1^TPA_1\)\dbE[\bar X]\\
\ns\ds\q=-K_0^{-1}(B^TP+B_1^TPA_1)\bar X.\ea$$
This recovers the result for classical LQ problem (\cite{Yong-Zhou
1999}).

\ms

\section{A Modification of Standard LQ Problems.}

In this section, we are going to look at a special case which was
mentioned in the introduction. For convenience, let us rewrite the
state equation here:
\bel{}\left\{\ba{ll}
\ns\ds
dX(s)=\[A(s)X(s)+B(s)u(s)\]ds+\[A_1(s)X(s)+B_1(s)u(s)\]dW(s),\\
\ns\ds X(0)=x,\ea\right.\ee
with the cost functional:
\bel{}J_0(x;u(\cd))=\dbE\[\int_0^T\(\lan Q_0(s)X(s),X(s)\ran+\lan
R_0(s)u(s),u(s)\ran\)ds+\lan G_0X(T),X(T)\ran\].\ee
Classical LQ problem can be stated as follows.

\ms

\bf Problem (LQ). \rm For any given $x\in\dbR^n$, find a $\bar
u(\cd)\in\cU[0,T]$ such that
\bel{}J_0(x;\bar u(\cd))=\inf_{u(\cd)\in\cU[0,T]}J_0(x;u(\cd)).\ee

\ms

The following result is standard (see \cite{Yong-Zhou 1999}).

\ms

\bf Theorem 5.1. \sl Let {\rm(H1)} hold and
\bel{}Q_0(s)\ge0,\q R_0(s)\ge\d I,\qq s\in[0,T];\qq G_0\ge0.\ee
Then Problem {\rm(LQ)} admits a unique optimal pair $(\bar
X_0(\cd),\bar u_0(\cd))$. Moreover, the following holds:
\bel{feedback1}\bar
u_0(s)=-\[R_0(s)+B_1(s)^TP_0(s)B(s)\]^{-1}\[B(s)^TP_0(s)+B_1(s)^TP_0(s)A_1(s)\]\bar
X_0(s),\q s\in[0,T],\ee
where $P_0(\cd)$ is the solution to the following Riccati equation:
\bel{Riccati P0}\left\{\2n\ba{ll}
\ns\ds
P_0'+P_0A+A^TP_0+A_1^TP_0A_1+Q_0\\
\ns\ds\qq-(P_0B+A_1^TP_0B_1)(R_0+B_1^TP_0B_1)^{-1}
(B^TP_0+B_1^TP_0A_1)=0,\q s\in[0,T],\\
\ns\ds P_0(T)=G_0.\ea\right.\ee
and $\bar X(\cd)$ is the solution to the following closed-loop
system:
\bel{}\left\{\3n\ba{ll}
\ns\ds d\bar
X_0(s)=\[A-B(R_0+B_1^TP_0B_1)^{-1}(B^TP_0+B_1^TP_0A)\]\bar
X(s)ds\\
\ns\ds\qq\qq+\[A_1-B_1(R_0+B_1^TP_0B_1)(B^TP_0+B_1^TP_0A_1)\]\bar X(s)dW(s),\q s\in[0,T],\\
\ns\ds\bar X_0(0)=x.\ea\right.\ee

\ms

\rm

We now introduce the following modified cost functional:
\bel{}\ba{ll}
\ns\ds\h J_0(x;u(\cd))=\dbE\[\int_0^T\(\lan Q_0(s)X(s),X(s)\ran+\lan
R_0(s)u(s),u(s)\ran\)ds+\lan
G_0X(T),X(T)\ran\]\\
\ns\ds\qq\qq+\dbE\[\int_0^T\(q(s)\var[X(s)]+\rho(s)\var[u(s)]\)ds
+g\var[X(T)\]\\
\ns\ds=\dbE\Big\{\int_0^T\[\lan\(Q_0(s)+q(s)I\)X(s),X(s)\ran-q(s)\Big|\dbE[X(s)]\Big|^2+\lan
\(R_0(s)+\rho(s)I\)u(s),u(s)\ran\\
\ns\ds\qq\qq-\rho(s)\Big|\dbE[u(s)]\Big|^2\]ds+\lan
\(G_0+gI\)X(T),X(T)\ran-g\Big|\dbE[X(T)]\Big|^2\Big\},\ea\ee
with $q(\cd),\rho(\cd)\in L^\infty(0,T)$, $g\in[0,\infty)$ such that
\bel{}q(s),\rho(s)\ge0,\qq s\in[0,T].\ee
Also, of course, we assume that
\bel{}\int_0^T[q(s)+\rho(s)]ds+g>0.\ee
We want to compare the above Problem (LQ) with the following
problem:

\ms

\bf Problem (LQ)$'$. \rm For any given $x\in\dbR^n$, find a $\bar
u(\cd)\in\cU[0,T]$ such that
\bel{}\h J_0(x;\bar u(\cd))=\inf_{u(\cd)\in\cU[0,T]}\h
J_0(x;u(\cd)).\ee
We refer to the above Problem (LQ)$'$ as a {\it modified} LQ
problem. This is a special case of Problem (MF) with
$$\left\{\ba{ll}
\ns\ds\h A=\h A_1=0,\q\h B=\h B_1=0,\\
\ns\ds Q=Q_0+qI,\q\h Q=-qI,\q R=R_0+\rho I,\q\h R=-\rho I,\q
G=G_0+gI,\q\h G=-gI.\ea\right.$$
Then the Riccati equations are
\bel{Riccati P1}\left\{\2n\ba{ll}
\ns\ds P'+PA+A^TP+A_1^TPA_1+Q_0+qI\\
\ns\ds\qq\qq-(PB+A_1^TPB_1)(R_0+\rho I+B_1^TPB_1)^{-1}
(B^TP+B_1^TPA_1)=0,\qq s\in[0,T],\\
\ns\ds P(T)=G_0+gI,\ea\right.\ee
and
%
%\bel{Riccati Pi1}\left\{\3n\ba{ll}
%
%\Pi'+\Pi\[A-B(R_0+B_1^TPB_1)^{-1}B_1^TPA_1\]+\[A^T-A_1^TPB_1(R_0+B_1^TPB_1)^{-1}B^T\]\Pi\\
%
%\ns\ds\q-\Pi B(R_0+B_1^TPB_1)^{-1}B^T\Pi+A_1^T\[P-PB_1(R_0+B_1^TPB_1)^{-1}B_1^TP\]A_1+Q_0=0,\q s\in[0,T],\\
%
%\ns\ds\Pi(T)=G_0.\ea\right.\ee
%
\bel{Riccati Pi1}\left\{\3n\ba{ll}
\ns\ds\Pi'\1n+\1n\Pi A\1n+\1n A^T\Pi\1n+\1n A_1^TPA_1\1n+\1n
Q_0\1n-\1n(\Pi B\1n+\1n A_1^TPB_1)(R_0\1n+\1n B_1^TPB_1)^{-1}
(B^T\Pi\1n+\1n B_1^TPA_1)\1n=\1n0,\q s\in[0,T],\\
\ns\ds\Pi(T)=G_0.\ea\right.\ee
The optimal control is given by
\bel{}\ba{ll}
\ns\ds\bar u=-(R_0+\rho I+B_1^TPB_1)^{-1}(B^TP+B_1^TPA_1)\(\bar
X-\dbE[\bar X]\)\\
\ns\ds\qq\qq-(R_0+B_1^TPB_1)^{-1}(B^T\Pi+B_1^TPA_1)\dbE[\bar
X],\ea\ee
and the closed-loop system reads
\bel{closed-loop}\left\{\ba{ll}
\ns\ds d\bar X=\Big\{\[A-B(R_0+\rho I+B_1^TPB_1)^{-1}(B^TP+B_1^TPA_1)\]\(\bar X-\dbE[\bar X]\)\\
\ns\ds\qq\qq+\[A-B(R_0+B_1^TPB_1)^{-1}\(B^T\Pi+B_1^TPA_1\)\]\dbE[\bar X]\Big\}ds\\
\ns\ds\qq\qq+\Big\{\[A_1-B_1K_0^{-1}(B^TP+B_1^TPA_1)\]\(\bar X-\dbE[\bar X]\)\\
\ns\ds\qq\qq+\[A_1-B_1K_1^{-1}\(B^T\Pi+B_1^TPA_1\)\]\dbE[\bar X]\Big\}dW(s),\\
\ns\ds\bar X(0)=x.\ea\right.\ee
By the optimality of $\bar u_0(\cd)$ and $\bar u(\cd)$, we have that
\bel{J<J}\lan P_0(0)x,x\ran=J_0(x;\bar u_0(\cd))\le J_0(x;\bar
u(\cd)),\ee
and
\bel{5.17}\ba{ll}
\ns\ds\lan\Pi(0)x,x\ran=J_0(x;\bar
u(\cd))+\dbE\[\int_0^T\(q(s)\var[\bar X(s)]+\rho(s)\var[\bar
u(s)]\)ds+g\var[\bar
X(T)]\]\\
\ns\ds=\h J_0(x;\bar u(\cd))\le\h J_0(x;\bar u_0(\cd))\\
\ns\ds=J_0(x;\bar
u_0(\cd))+\dbE\[\int_0^T\(q(s)\var[\bar X_0(s)]+\rho(s)\var[\bar u_0(s)]\)ds+g\var[\bar X_0(T)]\]\\
\ns\ds\le J_0(x;\bar u(\cd))+\dbE\[\int_0^T\(q(s)\var[\bar
X_0(s)]+\rho(s)\var[\bar u_0(s)]\)ds+g\var[\bar X_0(T)]\].\ea\ee
This implies
\bel{V<V}\ba{ll}
\ns\ds\dbE\[\int_0^T\(q(s)\var[\bar X(s)]+\rho(s)\var[\bar u(s)]\)ds+g\var[\bar X(T)]\]\\
\ns\ds\le\dbE\[\int_0^T\(q(s)\var[\bar X_0(s)]+\rho(s)\var[\bar
u_0(s)]\)ds+g\var[\bar X_0(T)]\].\ea\ee
Hence,
$$J_0(x;\bar u(\cd))-J_0(x;\bar u_0(\cd))$$
is the price for the decrease
$$\ba{ll}
\ns\ds\dbE\[\int_0^T\(q(s)\var[\bar X_0(s)]+\rho(s)\var[\bar u_0(s)]\)ds+g\var[\bar X_0(T)]\]\\
\ns\ds-\dbE\[\int_0^T\(q(s)\var[\bar X(s)]+\rho(s)\var[\bar
u(s)]\)ds+g\var[\bar X(T)]\ea$$
of the (weighted) variances of the optimal state-control pair $(\bar
X_0(\cd),\bar u_0(\cd))$. Moreover, (\ref{5.17}) further implies
that
$$\ba{ll}
\ns\ds J_0(x;\bar u(\cd))-J_0(x;\bar u_0(\cd))\\
\ns\ds\le\dbE\[\int_0^T\(q(s)\var[\bar X_0(s)]+\rho(s)\var[\bar u_0(s)]\)ds+g\var[\bar X_0(T)]\]\\
\ns\ds\qq-\dbE\[\int_0^T\(q(s)\var[\bar X(s)]+\rho(s)\var[\bar
u(s)]\)ds+g\var[\bar X(T)].\ea$$
The above roughly means that the amount increased in the cost is
``covered'' by the amount decreased in the weighted variance of the
optimal state-control pair.

\ms

We now look at a simple case to illustrate the above. Let us look at
a one-dimensional controlled linear SDE:
\bel{5.15}\left\{\ba{ll}
\ns\ds dX(s)=bu(s)ds+X(s)dW(s),\\
\ns\ds X(0)=x,\ea\right.\ee
with cost functionals:
\bel{5.16}J_0(x;u(\cd))=\dbE\[\int_0^T|u(s)|^2ds+g_0|X(T)|^2\],\ee
and
\bel{5.21}\ba{ll}
\ns\ds\h J_0(x;u(\cd))=\dbE\Big\{\int_0^T|u(s)|^2ds+g_0|X(T)|^2+g\var[X(T)]\Big\}\\
\ns\ds\qq\qq~=\dbE\Big\{\int_0^T|u(s)|^2ds+(g_0+g)|X(T)|^2-g\(\dbE[X(T)]\)^2\Big\},\ea\ee
where $g_0\ge0$ and $g>0$. As above, we refer to the optimal control
problem associated with (\ref{5.15}) and (\ref{5.16}) as the {\it
standard LQ problem}, and to that associated with (\ref{5.15}) and
(\ref{5.17}) as the {\it modified LQ problem}. The Riccati equation
for the standard LQ problem is
\bel{}\left\{\ba{ll}
\ns\ds
p_0'(s)+p_0(s)-b^2p_0(s)^2=0,\qq s\in[0,T],\\
\ns\ds p_0(T)=g_0.\ea\right.\ee
A straightforward calculation shows that
\bel{}p_0(s)={e^{T-s}g_0\over(e^{T-s}-1)b^2g_0+1}>0,\qq
s\in[0,T].\ee
The optimal control is
\bel{}\bar u_0(s)=-bp_0(s)\bar X_0(s),\qq s\in[0,T],\ee
and the closed-loop system is
\bel{}\left\{\ba{ll}
\ns\ds d\bar X_0(s)=-b^2p_0(s)\bar X_0(s)ds+\bar X_0(s)dW(s),\qq s\in[0,T],\\
\ns\ds\bar X_0(0)=x.\ea\right.\ee
Thus,
\bel{}\bar X_0(s)=xe^{-b^2\int_0^sp_0(\t)d\t-{1\over2}s+W(s)},\qq
s\in[0,T].\ee
Consequently,
\bel{}\dbE[\bar
X_0(T)]=xe^{-b^2\int_0^Tp_0(\t)d\t-{1\over2}T+{1\over2}T}=xe^{-b^2\int_0^Tp_0(\t)d\t},\ee
and
\bel{}\dbE[\bar
X_0(T)^2]=x^2e^{-2b^2\int_0^Tp_0(\t)d\t-T+2T}=x^2e^{-2b^2\int_0^Tp_0(\t)d\t+T}.\ee
Hence,
\bel{}\var[\bar X_0(T)]=\dbE[\bar X_0(T)^2]-\(\dbE[\bar
X_0(T)]\)^2=x^2e^{-2b^2\int_0^Tp_0(\t)d\t}(e^T-1).\ee
Also, the optimal expected cost is
\bel{}J_0(x;\bar
u_0(\cd))=p_0(0)x^2={e^Tg_0\over(e^T-1)b^2g_0+1}x^2.\ee
Next, for the modified LQ peoblem, the Riccati equations are:
\bel{}\left\{\ba{ll}
\ns\ds
p'(s)+p(s)-b^2p(s)^2=0,\qq s\in[0,T],\\
\ns\ds p(T)=g_0+g,\ea\right.\ee
and
\bel{}\left\{\ba{ll}
\ns\ds\pi'(s)-b^2\pi^2(s)+p(s)=0,\qq s\in[0,T],\\
\ns\ds\pi(T)=g_0.\ea\right.\ee
Clearly,
\bel{}p(s)={e^{T-s}(g_0+g)\over(e^{T-s}-1)b^2(g_0+g)+1}>{e^{T-s}g_0\over(e^{T-s}-1)b^2g_0+1}=p_0(s)>0,\qq
s\in[0,T].\ee
We now show that
\bel{p0<pi<p}p_0(s)<\pi(s)<p(s),\qq s\in[0,T].\ee
In fact,
$$\left\{\ba{ll}
\ns\ds{d\over
ds}[\pi(s)-p_0(s)]-b^2[\pi(s)+p_0(s)][\pi(s)-p_0(s)]+p(s)-p_0(s)=0,\qq
s\in[0,T],\\
\ns\ds\pi(T)-p_0(T)=0.\ea\right.$$
Thus,
$$\pi(s)-p_0(s)=\int_s^Te^{-\int_s^tb^2[\pi(\t)+p_0(\t)]d\t}[p(t)-p_0(t)]dt>0,\qq s\in[0,T).$$
Next,
$$\left\{\ba{ll}
\ns\ds{d\over ds}[p(s)-\pi(s)]-b^2[p(s)+\pi(s)][p(s)-\pi(s)]=0,\qq
s\in[0,T],\\
\ns\ds p(T)-\pi(T)=g.\ea\right.$$
Hence,
$$p(s)-\pi(s)=e^{-\int_s^Tb^2[p(\t)+\pi(\t)]d\t}g>0,\qq s\in[0,T].$$
This proves (\ref{p0<pi<p}). Note that the optimal control of
modified LQ problem is given by
$$\bar u(s)=-bp(s)\(\bar X(s)-\dbE[\bar
X(s)]\)-b\pi(s)\dbE[\bar X(s)],\qq s\in[0,T],$$
and the closed-loop system is
\bel{}\left\{\ba{ll}
\ns\ds d\bar X(s)=-\[b^2p(s)\(\bar X(s)-\dbE[\bar
X(s)]\)+b^2\pi(s)\dbE[\bar X(s)]\]ds+\bar X(s)dW(s),\qq
s\in[0,T],\\
\ns\ds\bar X(0)=x.\ea\right.\ee
Thus,
$${d\over ds}\(\dbE[\bar X(s)]\)=-b^2\pi(s)\dbE[\bar X(s)],$$
which leads to
\bel{}\dbE[\bar X(s)]=e^{-b^2\int_0^s\pi(\t)d\t}x,\qq s\in[0,T].\ee
On the other hand, by It\^o's formula,
$$d[\bar X^2]=\Big\{-2\[b^2p\bar X\(\bar
X-\dbE[\bar X]\)+b^2\pi\bar X\dbE[\bar X]\]+\bar
X^2\Big\}ds+[\cds]dW.$$
Then
$$\ba{ll}
\ns\ds d\(\dbE[\bar X^2]\)=\Big\{-2b^2p\[\dbE[\bar X^2]-\(\dbE[\bar
X]\)^2\]-2b^2\pi\(\dbE[\bar X]\)^2+\dbE[\bar
X^2]\Big\}ds\\
\ns\ds\qq=\Big\{(1-2b^2p)\dbE[\bar X^2]+2b^2(p-\pi)\(\dbE[\bar
X]\)^2\Big\}ds\\
\ns\ds\qq=\Big\{(1-2b^2p)\dbE[\bar
X^2]+2b^2(p-\pi)e^{-2b^2\int_0^s\pi(\t)d\t}x^2\Big\}ds.\ea$$
Hence,
$$\ba{ll}
\ns\ds\dbE[\bar
X^2(s)]=e^{s-2b^2\int_0^sp(\t)d\t}x^2\[1+\int_0^se^{-t+2\int_0^tp(\t)d\t}2b^2[p(t)-\pi(t)]
e^{-2b^2\int_0^t\pi(\t)d\t}dt\]\\
\ns\ds=e^{s-2b^2\int_0^sp(\t)d\t}x^2\[1+\int_0^se^{-t+2b^2\int_0^t[p(\t)-\pi(\t)]d\t}2b^2[p(t)-\pi(t)]
dt\]\\
%
%\ns\ds=e^{s-2b^2\int_0^sp(\t)d\t}x^2\[1+\int_0^se^{-t}d\(e^{2b^2\int_0^t[p(\t)-(1+\hat r_0)\pi(\t)]d\t}\)\]\\
%
\ns\ds=e^{s-2b^2\int_0^sp(\t)d\t}x^2\[e^{-s+2b^2\int_0^s[p(\t)-\pi(\t)]d\t}+\int_0^se^{-t+2b^2
\int_0^t[p(\t)-\pi(\t)]d\t}dt\]\\
\ns\ds=\[e^{-2b^2\int_0^s\pi(\t)d\t}+\int_0^se^{s-t-2b^2\int_t^sp(\t)d\t-2b^2\int_0^t\pi(\t)d\t}dt\]x^2.\ea$$
It follows that
$$\ba{ll}
\ns\ds\var[\bar
X(s)]=\[e^{-2b^2\int_0^s\pi(\t)d\t}+\int_0^se^{s-t-2b^2\int_t^sp(\t)d\t-2b^2\int_0^s\pi(\t)d\t}dt\]x^2-\(e^{-2b^2\int_0^t\pi(\t)d\t}\)x^2\\
\ns\ds\qq\qq\;=\[\int_0^se^{s-t-2b^2\int_t^sp(\t)d\t-2b^2\int_0^t\pi(\t)d\t}dt\]x^2.\ea$$
Consequently, noting $p_0(s)<p(s)$, we have
\bel{5.37}\ba{ll}
\ns\ds\var[\bar
X(T)]\1n=\1n\[\1n\int_0^T\2ne^{T-t-2b^2\int_t^Tp(\t)d\t-2b^2\int_0^tp_0(\t)d\t}dt\]x^2
\1n<\1n\[\1n\int_0^T\2ne^{T-t-2b^2\int_t^Tp_0(\t)d\t-2b^2\int_0^tp_0(\t)d\t}dt\]x^2\\
\ns\ds\qq\qq~=\[\int_0^Te^{T-t-2b^2\int_0^Tp_0(\t)d\t}dt\]x^2=\[e^{-2b^2\int_0^Tp_0(\t)d\t}
\int_0^Te^{T-t}dt\]x^2=\var[\bar X_0(T)].\ea\ee
On the other hand, we claim that
\bel{wt p>p}p_0(s)>{g_0\over g_0+g}p(s),\qq s\in[0,T).\ee
In fact, by letting $\wt p={g_0\over g_0+g}p$, we have
$$\left\{\ba{ll}
\ns\ds\wt p'(s)+\wt p(s)-b^2\wt p(s)^2-{g_0g\over(g_0+g)^2}p(s)^2=0,\\
\ns\ds\wt p(T)=g_0.\ea\right.$$
Then
$$\left\{\ba{ll}
\ns\ds[p_0'(s)-\wt p'(s)]+[p_0(s)-\wt p(s)]-b^2[p_0(s)+\wt
p(s)][p_0(s)-\wt p(s)]+{g_0g\over(g_0+g)^2}p(s)^2=0,\\
\ns\ds p_0(T)-\wt p(T)=0.\ea\right.$$
This leads to (\ref{wt p>p}). Consequently,
$$\ba{ll}
\ns\ds J_0(x;\bar u(\cd))+\var[\bar X(T)]=\h J_0(x;\bar
u(\cd))=\pi(0)x^2\le p(0)x^2\le{g_0+g\over g_0}p_0(0)x^2={g_0+g\over
g_0}J(x;\bar u_0(\cd)).\ea$$
Therefore,
$$0\le J_0(x;\bar u(\cd))-J_0(x;\bar u_0(\cd))+g\var[\bar X(T)]\le{g\over
g_0}J(x;\bar u_0(\cd)).$$
The above gives an upper bound for the cost increase in order to
have a smaller $\var[X(T)]$. Taking into account of (\ref{5.37}), we
see that it is a very good trade-off to consider the modified LQ
problem if one wishes to have a smaller $\var[X(T)]$. It is possible
to more carefully calculate the price difference
$J_0(x;u(\cd))-J_0(x;\bar u_0(\cd))$. We omit the details here.

\ms

Also, it is possible to calculate the situation of including
$\var[u(s)]$ and/or $\var[X(s)]$ in the integrand of the modified
cost functional. The details are omitted here as well.

\ms

To conclude this paper, let us make some remarks. First of all, we
have presented some results on the LQ problem for MF-SDEs with
deterministic coefficients. Optimal control is represented by a
state feedback form involving both $X(\cd)$ and $\dbE[X(\cd)]$, via
the solutions of two Riccati equations. Apparently, there are many
problems left unsolved (and some of them might be challenging). To
mention a few, one may consider the case of infinite-horizon
problem, following the idea of \cite{Wu-Zhou 2001}, and more
interestingly, the case of random coefficients (for which one might
have to introduce some other techniques since the approach used in
this paper will not work). We will continue our study and report new
results in our future publications.

\ms

\end{document}